\documentclass[10pt,a4paper]{article}

\usepackage{mystyle}
\usepackage{mythmstyle_report}

\begin{document}

\tikzset{->-/.style={decoration={
markings,
mark=at position #1 with {\arrow{>}}},postaction={decorate}}}

\title{Cayley graphs and symmetric interconnection networks}

\author{Ashwin~Ganesan%
  \thanks{These lecture notes were published as: A. Ganesan, ``Cayley graphs and symmetric interconnection networks,'' Proceedings of the Pre-Conference Workshop on Algebraic and Applied Combinatorics (PCWAAC 2016), 31st Annual Conference of the Ramanujan Mathematical Society, pp. 118--170, Trichy, Tamilnadu, India, June 2016.
  Correspondence address: 53 Deonar House, Deonar Village Road, Deonar, Mumbai 400088, Maharashtra, India. Email:  
\texttt{ashwin.ganesan@gmail.com}.}
}

\date{}
\vspace{10cm}

\maketitle

%
\begin{abstract}
\noindent  
These lecture notes are on automorphism groups of Cayley graphs and their applications to optimal fault-tolerance of some interconnection networks.  We first give an introduction to automorphisms of graphs and an introduction to Cayley graphs.   We then discuss automorphism groups of Cayley graphs.  We prove that the vertex-connectivity of edge-transitive graphs is maximum possible. We investigate the automorphism group and vertex-connectivity of some families of Cayley graphs that have been considered for interconnection networks; we focus on the hypercubes, folded hypercubes, Cayley graphs generated by transpositions, and Cayley graphs from linear codes.  New questions and open problems are also discussed.
\end{abstract}


\vskip 0.3in

\tableofcontents

\vspace{+0.5cm}

\section{Introduction}

In these lecture notes, we discuss the automorphism groups and fault-tolerance of some families of Cayley graphs.  These notes are organized as follows. We first introduce the topic of automorphisms of graphs in Section~\ref{sec:aut:graphs} and Cayley graphs in  Section~\ref{sec:intro:cayley:graphs}.  In Section~\ref{sec:intro:icn} we discuss interconnection networks of supercomputers; in particular, we discuss the degree-diameter problem and low-dimensional networks such as torus graphs which are used as the interconnection network of some of the fastest supercomputers. 

In Section~\ref{sec:aut:cayley}, we study the automorphism groups of Cayley graphs. We prove that all Cayley graphs are vertex-transitive and that the converse is not true in general - the Petersen graph is a vertex-transitive graph which is not a Cayley graph. We prove Sabidussi's theorem that a graph is a Cayley graph if and only if its automorphism group contains a regular subgroup. We define normal Cayley graphs and mention some recent results.  

If a graph possesses certain symmetries, then the graph has high connectivity properties that make it a good topology to use for interconnection networks.  In Section~\ref{sec:connec:of:transitive}, we mention the fundamental results on connectivity of transitive graphs. We first define the notions of vertex-connectivity, fault-tolerance and edge-connectivity of graphs.  We prove Whitney's theorem that the vertex-connectivity of a graph is at most its edge-connectivity.  A graph is said to have optimal fault-tolerance (or vertex-connectivity that is maximum possible) if its vertex-connectivity is equal to its minimum degree.  We use Watkins' theory of atomic parts to prove that an edge-transitive graph has vertex-connectivity that is maximum possible, and we state some other known results on connectivity of symmetric graphs.

In Section~\ref{sec:kappaeqdelta}, we prove that the vertex-connectivity of many families of Cayley graphs is maximum possible.  We focus on the hypercube graphs, the folded hypercubes, Cayley graphs generated by transpositions, and Cayley graphs from linear codes.  For the first two families of graphs, we give several proofs of their optimal fault-tolerance.

Throughout, $X=(V,E)$ denotes a simple, undirected graph with vertex set $V=V(X)$ and edge set $E=E(X)$.  If $H$ is a group and $S \subseteq H$, then $\Cay(H,S)$ denotes the Cayley graph of $H$ with respect to $S$, and its definition is given below.  All graphs and groups studied here are finite.    We assume the reader is familiar with basic notions from algebra (cf. Gallian \cite{Gallian:1990}) and graph theory (cf. Bollob\'as \cite{Bollobas:1998}).   The minimum degree of a vertex of $X$ is denoted $\delta(X)$, and a $k$-regular graph is said to have valency $k$.  For a vertex $v \in V$, the $i$th layer of $X$ with respect to $v$, denoted by $X_i(v)$, is the set of vertices of $X$ whose distance to the vertex $v$ is exactly $i$.  Thus, $X_0(v) = \{v\}$ and $X_1(v)$ is the set of neighbors of $v$.

\subsection{Automorphisms of graphs}
\label{sec:aut:graphs}

In this section, we define the automorphism group of a graph and give some examples.  We first recall some basic definitions and terminology on group actions. Standard references on permutation groups are Cameron \cite{Cameron:1999}, Dixon and Mortimer \cite{Dixon:Mortimer:1996}, and Wielandt \cite{Wielandt:1964}. 

Given a finite set $V$, the set of all bijections from $V$ to itself, denoted by $\Sym(V)$, forms a permutation group called the full symmetric group on $V$.  If $V = \{1,\ldots,n\}$, then $\Sym(V)$ is abbreviated $S_n$.   Let $G$ be a finite group and let $V$ be a finite set.  Let $1$ denote the identity element of the group $G$.  Suppose there is a map from $V \times G$ into $V$ which satisfies the following two conditions: (a) $u^1=u$, for all $u \in V$, and (b) $(u^g)^h = u^{(gh)}$, for all $u \in V$ and $g, h \in G$.  Then, we say this map defines an {\em action} of $G$ on $V$.  Every action of $G$ on $V$ induces a homomorphism from $G$ into $\Sym(V)$, and conversely. If $g \in \Sym(V)$, then there is a corresponding action of $g$ on the set $V^{\{2\}}$ of all 2-subsets of $V$, defined by the rule $\{u,v\}^g := \{u^g, v^g\}$ for each $\{u,v\} \in V^{\{2\}}$ and each $g \in G$.   Similarly, if a group $G$ acts on $V$, then $G$ acts naturally on the cartesian product $V \times V$. This induced action is defined by $(u,v)^g := (u^g, v^g)$ for $(u,v) \in V \times V$ and $g \in G$.  

Let $V$ be a finite set and let $G$ be a finite group acting on $V$.  The orbit of $G$ containing $v \in V$ is the set $v^G := \{v^g: v \in V\}$ $u \in V$.  The action of $G$ on $V$ partitions $V$ into orbits.  We say that $G$ acts {\em  transitively} on $V$ if for all $u,v \in V$, there is an element $g \in G$ such that $u^g = v$.  Equivalently, $G$ is transitive on $V$ if the action has a single orbit.  The group $G$ acts {\em regularly} on $V$ if for all $u, v \in V$, there is a unique element $g \in G$ such that $u^g = v$.  The stabilizer in $G$ of the point $v$, denoted by $G_v$, is the set of elements in $G$ that fix $v$, i.e.,  $G_v := \{g \in G: v^g = v\}$.  A fundamental result (often referred to as the Orbit-Stabilizer Lemma) in the theory of permutation groups states that  for each $v \in V$, the length of the orbit of $G$ containing $v$ is equal to the index of the stabilizer in $G$ of $v$, i.e., that $|G_v|~|v^G| = |G|$. 

Let $X=(V,E)$ be a simple, undirected graph with vertex set $V=V(X)$ and edge set $E=E(X)$.  An {\em automorphism} of the graph $X$ is a permutation of the vertex set $V(X)$ which preserves adjacency. In other words, a permutation $g \in \Sym(V)$ is an automorphism of $X$ if $\{u,v\} \in E$ iff $\{u^g, v^g\} \in E$.  Thus, an automorphism of $X$ is an isomorphism from $X$ to itself.  The set of all automorphisms of the graph $X$, denoted by $\Aut(X)$, is called the {\em automorphism group of $X$}. In symbols, \[\Aut(X) := \{g \in \Sym(V): E^g = E\}.\]  While one can often obtain some automorphisms of a graph, it is often difficult to prove that one has obtained the (full) automorphism group.   Recall that the action of $g$ on $V$ naturally induces an action of $g$ on the set of all 2-subsets of $V$, and so $E^g$ is the image of the edge set $E$ under $g$.  A graph $X$ is said to be {\em vertex-transitive} if $\Aut(X)$ acts transitively on $V(X)$.  A graph $X$ is said to be {\em edge-transitive} if $\Aut(X)$ acts transitively on $E(X)$.

\begin{figure}
\begin{centering}
\begin{tikzpicture}[scale=3.5,auto]

\pgfmathsetmacro{\phi}{0.6}

\vertex[fill] (v12) at (0,1) [label=above:$12$] {};
\vertex[fill] (v34) at (0.95, 0.31) [label=right:$34$] {};
\vertex[fill] (v15) at (0.588, -0.81) [label=right:$15$] {};
\vertex[fill] (v23) at (-0.588,-0.81) [label=below:$23$] {};
\vertex[fill] (v45) at (-0.95, 0.31) [label=left:$45$] {};
\draw (v12) -- (v34) -- (v15) -- (v23) -- (v45) -- (v12);

\vertex[fill] (v35) at (0 * \phi,1 * \phi) [label=right:$35$] {};
\vertex[fill] (v25) at (0.95 * \phi, 0.31 * \phi) [label=below:$25$] {};
\vertex[fill] (v24) at (0.588 * \phi, -0.81 * \phi) [label=right:$24$] {};
\vertex[fill] (v14) at (-0.588 * \phi,-0.81 * \phi) [label=below:$14$] {};
\vertex[fill] (v13) at (-0.95 * \phi, 0.31 * \phi) [label=below:$13$] {};
\draw (v12) -- (v34) -- (v15) -- (v23) -- (v45) -- (v12);
\draw (v35) -- (v24) -- (v13) -- (v25) -- (v14) -- (v35);
\draw (v12) -- (v35);
\draw (v34) -- (v25);
\draw (v15) -- (v24);
\draw (v23) -- (v14);
\draw (v45) -- (v13);

\end{tikzpicture}
\caption{The Petersen graph}
\label{fig:Petersen:graph}
\end{centering}
\end{figure}
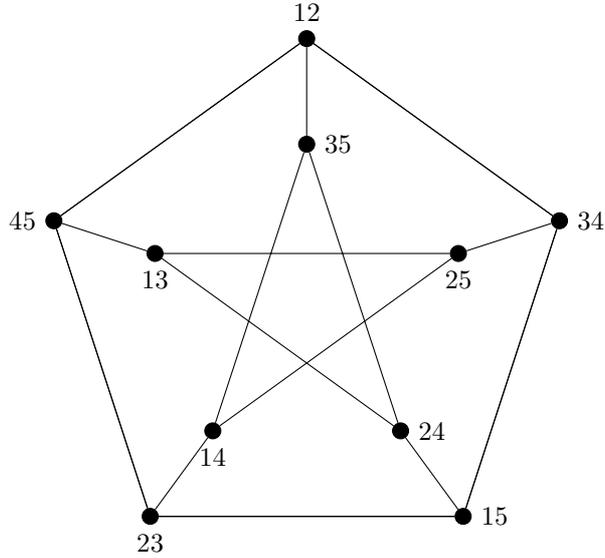

The Petersen graph $X$ is defined as follows.  The vertex set $V$ is the set of all 2-element subsets of the set $[5]=\{1,2,3,4,5\}$. Two vertices $A$ and $B$ are adjacent in this graph whenever $A \cap B = \phi$; see Figure~\ref{fig:Petersen:graph}, where vertex $\{i,j\}$ is labeled $ij$ for brevity.  It is clear that $S_5$ acts naturally on the 2-subsets of $[5]$ as a group of automorphisms of $X$, for if $g \in S_5$ and $A, B \in V$, then $|A \cap B| = |A^g \cap B^g|$.  Thus, $S_5$ is isomorphic to a subgroup of the automorphism group of the Petersen graph.  It turns out that $S_5$ is the full automorphism group of the Petersen graph. The nontrivial part of this result is to show that the Petersen graph has no other automorphisms besides those induced by $S_5$. We give two proofs of this result. 

\begin{Theorem} \label{thm:autgroup:Petersen}
 The automorphism group of the Petersen graph is isomorphic to $S_5$. 
\end{Theorem}

\noindent \emph{First proof:} Let $X$ be the Petersen graph, shown in Figure~\ref{fig:Petersen:graph}. As explained above, $S_5$ acts on $V(X)$ as a group of automorphisms of $X$.  Thus, $S_5$ is isomorphic to a subgroup of $\Aut(X)$.  We need to show that $X$ has no other automorphisms besides those induced by $S_5$.  Let $g$ be an automorphism of $X$.  It suffices to show that there exist elements $w, x, y, z \in S_5$ such that the action of $g$ on $V(X)$ composed with the action of $wxyz$ on $V(X)$ gives the identity permutation of $V(X)$.  

Suppose the automorphism $g$ takes vertex $\{1,2\}$ to vertex $\{i,j\}$. Then take $w$ to be any permutation in $S_5$ that maps $i$ to 1 and $j$ to 2.  Observe that $gw$ is an automorphism of $X$ that fixes vertex $\{1,2\}$. Consider the action of $gw$ on the three vertices $\alpha = \{3,4\}$, $\beta = \{3,5\}$, $\gamma = \{4,5\}$. If $gw$ takes $\alpha$ to $\beta$ or $\gamma$, then there exists an $x \in S_5$ such that $gwx$ fixes both $\{1,2\}$ and $\alpha$. For example, if $gw$ takes $\alpha=\{3,4\}$ to $\beta=\{3,5\}$, then take $x=(1)(2)(3)(4,5)$; if $gw$ fixes $\alpha$, then take $x=1$. Continuing in this manner, we see that there exists a $y \in S_5$ such that $gwxy$ fixes each of the vertices $\{1,2\}, \alpha, \beta$ and $\gamma$.  If $gwxy$ is the identity permutation of $V(X)$, then we are done.  Otherwise, since $gwxy$ is an automorphism of $X$ and preserves adjacencies, it can be seen that $gwxy$ must interchange vertices $\{1,i\}$ and $\{2,i\}$ for each $i=3,4,5$.  Taking $z=(1,2) \in S_5$ , we get that $gwxyz$ is the identity permutation of $V(X)$.  
\qed

Given a graph $X=(V,E)$, the {\em line graph} of $X$, denoted by $L(X)$, is defined to be the graph with vertex set $E$, and two vertices $e, f \in E$ are adjacent in $L(X)$ whenever $e, f$ are incident in $X$.  Thus, the Petersen graph is the complement of the line graph of $K_5$.  Each star $K_{1,r}$ in a graph $X$ gives rise to a clique $K_r$ in $L(X)$.  Conversely, every clique in $L(X)$ of size 4 or more is induced by a star in $X$.  

\bigskip \noindent \emph{Second proof of Theorem~\ref{thm:autgroup:Petersen}:} 
The Petersen graph is isomorphic to the complement of $L(K_5)$.  Since a graph and its complement have the same automorphism group, it suffices to find the automorphism group of the graph $L(K_5)$.  Let $X := L(K_5)$ and let $G:=\Aut(X)$.  The automorphism group of the complete graph $K_5$ is $S_5$.  An automorphism of a graph induces an automorphism of its line graph, and so $S_5$ also acts as a group of automorphisms of $L(K_5)$.  Thus, $S_5$ is isomorphic to a subgroup of $G$.  

To prove that $G \cong S_5$, we consider a set of 5 substructures in the graph $X$ for the domain of the action of its automorphism group $G$.   The graph $X=L(K_5)$ has exactly five maximum cliques $\mathcal{C}_1,\ldots,\mathcal{C}_5$, each clique  corresponding to a star $K_{1,4}$ in $K_5$.  An automorphism of a graph must permute its maximum cliques among themselves.  Consider the action of $G$ on the set $\{\mathcal{C}_1,\ldots,\mathcal{C}_5\}$ of maximum cliques of $X$.   To prove $|G| \le 120$, it suffices to show that the kernel of this action is trivial.  In the graph $X$, each vertex lies in exactly two maximum cliques, and any two maximum cliques intersect in exactly one point.  Thus, if $g \in G$ takes each maximum clique to itself, then $g$ also fixes the unique common point of any two cliques, and hence $g$ fixes each vertex of $X$. Thus, the kernel of the action is trivial. This completes the proof. 
\qed

\begin{Exercise}
 Let $G$ be the automorphism group of the Petersen graph $X$.  Show that the stabilizer in $G$ of a vertex has 3 orbits, of lengths 1, 3 and 6, respectively.  
\end{Exercise}

\begin{Example}
The complete bipartite graph $K_{3,3}$ is both vertex-transitive and edge-transitive.  The graph $K_{3,4}$ is edge-transitive but not vertex-transitive.  The star graph $K_{1,5}$ is also edge-transitive but not vertex-transitive.  
\qed
\end{Example}

\begin{Example}
The pentagonal prism $L_5$ (cf. Figure~\ref{fig:pentagonal:prism}) is vertex-transitive but not edge-transitive since it has two types of edges: those that are a part of a 5-cycle and those that are not.  The action of the automorphism group $\Aut(L_5)$ on the edge set $E(L_5)$ has two orbits.
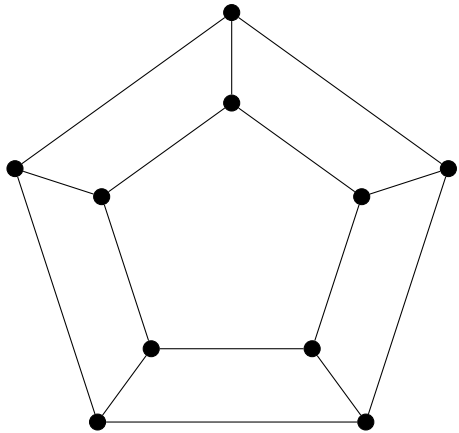
\begin{figure}
\begin{centering}
\begin{tikzpicture}[scale=3,auto]

\pgfmathsetmacro{\phi}{0.6}

\vertex[fill] (x1) at (0,1) {};
\vertex[fill] (x2) at (0.95, 0.31) {};
\vertex[fill] (x3) at (0.588, -0.81) {};
\vertex[fill] (x4) at (-0.588,-0.81) {};
\vertex[fill] (x5) at (-0.95, 0.31) {};
\draw (x1) -- (x2) -- (x3) -- (x4) -- (x5) -- (x1);

\vertex[fill] (y1) at (0 * \phi,1 * \phi) {};
\vertex[fill] (y2) at (0.95 * \phi, 0.31 * \phi) {};
\vertex[fill] (y3) at (0.588 * \phi, -0.81 * \phi) {};
\vertex[fill] (y4) at (-0.588 * \phi,-0.81 * \phi) {};
\vertex[fill] (y5) at (-0.95 * \phi, 0.31 * \phi) {};
\draw (y1) -- (y2) -- (y3) -- (y4) -- (y5) -- (y1);
\draw (x1) -- (y1);
\draw (x2) -- (y2);
\draw (x3) -- (y3);
\draw (x4) -- (y4);
\draw (x5) -- (y5);

\end{tikzpicture}
\caption{The pentagonal prism}
\label{fig:pentagonal:prism}
\end{centering}
\end{figure}
\qed
\end{Example}

Let $H$ and $K$ be subgroups of $G$.  Recall that $G$ is said to be the {\em direct product of $H$ and $K$}, denoted by $G = H \times K$, if $H$ and $K$ are both normal subgroups of $G$, $H \cap K = 1$ and $G=HK$.  The group $G$ is said to be the {\em semidirect product of $H$ by $K$}, denoted by $G = H \rtimes K$, if $H$ is a normal subgroup of $G$, $H \cap K =1$ and $G=HK$.

\begin{Exercise}
Prove that the automorphism group of the pentagonal prism $L_5$ is the dihedral group $D_{20} = D_{10} \times C_2$ of order 20.
 \qed
\end{Exercise}

One can also define the automorphism group of more general combinatorial structures.  For example, if $V$ is a finite set and $\mathcal{R}$ is a subset of $V \times \cdots \times V$ ($k$ times), then the automorphism group of the relational structure $(V,\mathcal{R})$ is defined to be the set of all permutations of $V$ that takes $\mathcal{R}$ to itself: $\Aut(V, \mathcal{R}) := \{ g \in \Sym(V): \mathcal{R}^g = \mathcal{R} \}$.  In particular, when $k=2$, we get the definition for the automorphism group of a digraph with vertex set $V$ and arc set $\mathcal{R} \subseteq V \times V$. The automorphism group of an edge-colored digraph $D=(V,\mathcal{A})$ can be defined similarly when $V$ is the vertex set of the digraph and $\mathcal{A}$ is a set of 3-tuples of the form $(u,v,c)$, where $c$ is the color of the arc $(u,v)$. 

\begin{Example}
 Let $X$ be the complete bipartite graph $K_{3,3}$.  The complement graph $\overline{X}$ consists of the disjoint union of two 3-cliques.  An automorphism of $\overline{X}$ can permute the vertices of the first clique among themselves, can permute the vertices of the second clique among themselves, and can further possibly also interchange the two cliques.  Thus, $\Aut(X)$ is isomorphic to the semidirect product $(S_3 \times S_3) \rtimes C_2$ (cf. \cite[p. 46]{Dixon:Mortimer:1996}).  Observe that the automorphism group of $\overline{X}$ is isomorphic to the automorphism group of the structure $(\{1,\ldots,6\}, \{ \{1,2,3\},\{4,5,6\}\})$ which is a partition of a 6-element set into two equal-sized subsets.  
 \qed
\end{Example}

Let $X$ be a graph and let $L(X)$ denote the line graph of $X$. It is clear that if $X$ and $Y$ are isomorphic graphs, then $L(X)$ and $L(Y)$ are isomorphic.  The converse is not necessarily true: the graphs $K_{1,3}$ and $K_3$ are not isomorphic, but their line graphs are isomorphic (and equal to $K_3$).   Whitney \cite{Whitney:1932} proved that this pair is the only connected counterexample: if $X$ and $Y$ are connected graphs on 5 or more vertices such that $L(X) \cong L(Y)$, then $X \cong Y$.  In terms of automorphisms, it is clear that every automorphism of $X$ induces an automorphism of $L(X)$.  Whitney \cite{Whitney:1932} showed that we can go in the reverse direction: every automorphism of $L(X)$ is induced by a unique automorphism of $X$ if $X$ is a connected graph on 5 or more vertices.   

\begin{Theorem}
 Let $X$ be a connected graph on 5 or more vertices. Then the automorphism group of $X$ and the automorphism group of its line graph $L(X)$ are isomorphic. 
\end{Theorem}

A graph is said to be {\em asymmetric} if it has no nontrivial symmetries, i.e. if its automorphism group is trivial. Erd\"os and R\'enyi \cite{Erdos:Renyi:1963} proved that the proportion of graphs on $n$ vertices which have a trivial automorphism group approaches 1 as $n$ goes to infinity.

\begin{Theorem}
 Almost all graphs are asymmetric. 
\end{Theorem}

For further details on automorphism groups of graphs, we refer the reader to  
Lov\'asz\cite[Chapter 12]{Lovasz:2007}, Biggs \cite{Biggs:1993}, Godsil and Royle \cite{Godsil:Royle:2001}, Cameron \cite{Cameron:2005} and Babai \cite{Babai:1995}. 

\subsection{Cayley graphs}
\label{sec:intro:cayley:graphs}

Cayley graphs of groups are one way to generate many families of vertex-transitive graphs.  In this section, we introduce Cayley graphs through a few  examples.

\begin{Definition}
 Let $H$ be a finite group and let $S$ be a subset of $H$.  The {\em Cayley diagram (or Cayley digraph) of $H$ with respect to $S$}, denoted by $\Cay(H,S)$, is the digraph with vertex set $H$ and arc set $\{ (h, sh): h \in H, s \in S\}$.  
\end{Definition}

Throughout, $e$ denotes the identity element of the group $H$ and also the corresponding vertex of the digraph $\Cay(H,S)$. 

\begin{Example}
Consider the cyclic group $C_4 = \{e,a,a^2,a^3\}$, where $e$ denotes the identity group element. The Cayley diagram of $C_4$ with respect to the set $\{a\}$ is shown in Figure~\ref{fig:cay:digraph:c4}.  
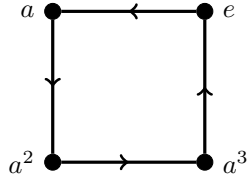
\begin{figure}
\begin{centering}
\begin{tikzpicture}[scale=1,auto]
\vertex[fill] (u1) at (2,2) [label=right:$e$] {};
\vertex[fill] (u2) at (0,2) [label=left:$a$] {};
\vertex[fill] (u3) at (0,0) [label=left:$a^2$] {};
\vertex[fill] (u4) at (2,0) [label=right:$a^3$] {};
\draw[very thick,->-=.5] (u1) to (u2);
\draw[very thick,->-=.5] (u2) to (u3);
\draw[very thick,->-=.5] (u3) to (u4);
\draw[very thick,->-=.5] (u4) to (u1);
\end{tikzpicture}
\caption{The Cayley digraph $\Cay(C_4 = \langle a \rangle, \{ a \})$.}
\label{fig:cay:digraph:c4}
\end{centering} 
\end{figure}
\qed
\end{Example}

\begin{Example}
The Cayley digraph of $S_3$ with respect to the set $S = \{a,b\}$, where $a=(123)$ and $b=(12)$, is shown in Figure~\ref{fig:cay:s3}.  For brevity, we will often write the cyclic permutation $(1,2,3)$ as $(123)$.  To construct the Cayley diagram, we start with the vertex $e$ corresponding to the identity group element.  By definition, the arcs of the Cayley digraph are precisely the pairs of the form $(h,sh)$ where $h \in S_3$ and $s \in S$.  Thus, the out-neighbors of the vertex $e$ are exactly the elements $a$ and $b$; see Figure~\ref{fig:cay:s3} for the arcs $(e,a)$ and $(e,b)$, respectively.  Next, we start at vertex $a=(12)$ and observe that its out-neighbors are $(123)a = (123)(12) = (23)$ and $(12)a = (12)(12) = e$.  We repeat this procedure of drawing two out-neighbors (corresponding to the two elements in $S$) from each new vertex.  Observe that when we start at $a^2$, we draw an arc from $a^2$ to two new vertices $ba^2$ and $a a^2=a^3$.  However, one of these new vertices, namely $a^3$, is equal to $e$. Hence, we {\em identify} the two vertices $a^3$ and $e$.  We repeat this process until each vertex has out-degree 2 and in-degree 2, and we identify vertices when possible.  This gives the Cayley digraph shown in Figure~\ref{fig:cay:s3}.  
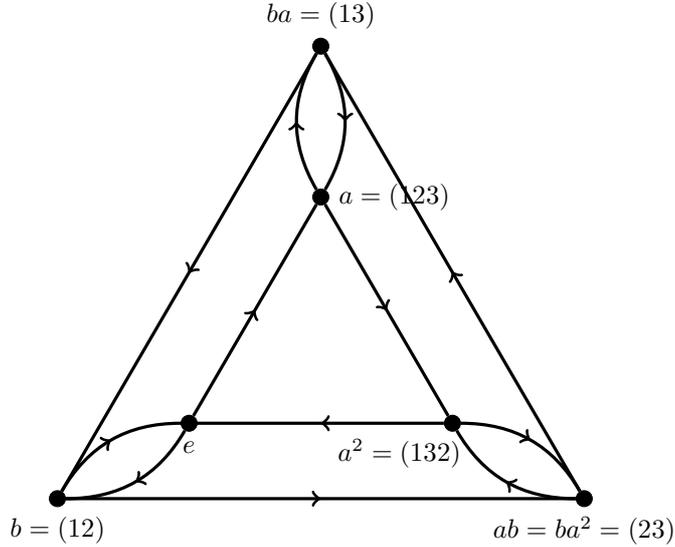
\begin{figure}
\begin{centering}
\begin{tikzpicture}[scale=2,auto]
\vertex[fill] (v1) at (0,1) [label=right:{$a=(123)$}] {};
\vertex[fill] (v2) at (0.866, -0.5) {};
\node at (0.866-0.35, -0.5-0.2) {$a^2=(132)$};
\vertex[fill] (v3) at (-0.866, -0.5) [label=below:{$e$}] {};
\vertex[fill] (x1) at (0,1*2) [label=above:{$ba=(13)$}] {};
\vertex[fill] (x2) at (0.866*2, -0.5*2) [label=below:{$ab=ba^2=(23)$}] {};
\vertex[fill] (x3) at (-0.866*2,-0.5*2) [label=below:{$b=(12)$}] {};

\draw[very thick, ->-=.5] (v1) to (v2);
\draw[very thick, ->-=.5] (v2) to (v3);
\draw[very thick, ->-=.5] (v3) to (v1);
\draw[very thick, ->-=.5] (x1) to (x3);
\draw[very thick, ->-=.5] (x3) to (x2);
\draw[very thick, ->-=.5] (x2) to (x1);

\draw[very thick, ->-=.5] (v3) to [bend left] (x3);
\draw[very thick, ->-=.5] (x3) to [bend left] (v3);

\draw[very thick, ->-=.5] (v2) to [bend left] (x2);
\draw[very thick, ->-=.5] (x2) to [bend left] (v2);

\draw[very thick, ->-=.5] (v1) to [bend left] (x1);
\draw[very thick, ->-=.5] (x1) to [bend left] (v1);

\end{tikzpicture}
\caption{The Cayley digraph $\Cay(S_3, \{(123), (12) \})$. }
\label{fig:cay:s3}
\end{centering}
\end{figure}

\qed
\end{Example}

We say $\Cay(H,S)$ is the Cayley digraph of $H$ generated by $S$.    Note that each vertex in the Cayley digraph has in-degree $|S|$ and out-degree $|S|$.   Each arc in the Cayley digraph is of the form $(h,sh)$ for some $h \in H, s \in S$.   The {\em edge-colored Cayley digraph} is obtained by assigning color (or edge label) $s$ to arc $(h,sh)$ in the Cayley digraph; in other words, it is the structure $(H, \mathcal{R})$, where $R = \{(h,sh,s): h \in H, s \in S \}$. 

If the identity group element $e$ belongs to $S$, then every vertex $h$ of the Cayley digraph will have a self-loop $(h,h)$.  Observe from Figure~\ref{fig:cay:s3} that since $b^2=e$ (i.e. $b$ is its own inverse), whenever there is an arc of color $b$ from vertex $u$ to $v$,  there is also an arc of color $b$ in the reverse direction (from $v$ to $u$).  Thus, if the generator set $S \subseteq H$ is closed under inverses, then $(u,v)$ is an arc of the Cayley digraph $\Cay(H,S)$ iff $(v,u)$ is an arc.  In this case, we replace the two arcs $(u,v)$ and $(v,u)$ in the digraph with a single undirected edge $\{u,v\}$.  Thus, when $e \notin S$ and $S=S^{-1}$, we consider $\Cay(H,S)$ to be a simple (no self-loops) undirected graph.   We summarize these observations with a formal definition: 

\begin{Definition}
 Let $H$ be a group and let $S$ be a subset of $H$ such that $e \notin S = S^{-1}$.  Then, the {\em Cayley graph $\Cay(H,S)$} is defined to be the graph with vertex set $H$ and edge set $\{ \{h, sh\}: h \in H, s \in S \}$.  
\end{Definition}

Observe that the Cayley graph $\Cay(H,S)$ is connected if and only if $S$ generates the group $H$.  Also, vertices $u, v \in H$ are adjacent in $\Cay(H,S)$ iff $uv^{-1} \in S$. 

\begin{Example}
Let $V_4$ be the Klein 4-group generated by the transpositions $a = (12)$ and $b = (34)$.  The Cayley graph $\Cay(V_4, \{a,b\})$ is shown in Figure~\ref{fig:cay:klein4}. 
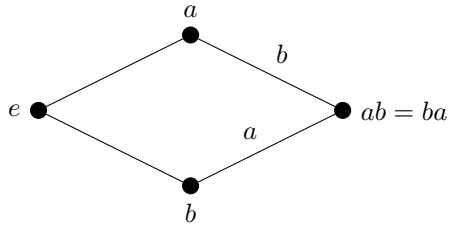
\begin{figure}
\begin{centering}
\begin{tikzpicture}[scale=1,auto]
\vertex[fill] (v1) at (6,1) [label=left:$e$] {};
\vertex[fill] (v2) at (8,0) [label=below:$b$] {};
\vertex[fill] (v3) at (8,2) [label=above:$a$] {};
\vertex[fill] (v4) at (10,1) [label=right:{$ab=ba$}] {};
\draw [-] (v1) edge (v2);
\draw [-] (v1) edge (v3);
\draw [-] (v2) edge node {$a$} (v4);
\draw [-] (v3) edge node {$b$} (v4);
\end{tikzpicture}
\caption{The Cayley graph of the Klein 4-group}
\label{fig:cay:klein4}
\end{centering} 
\end{figure}
\qed
\end{Example}

The {\em $n$-dimensional hypercube graph $Q_n$} is defined to be the graph with vertex set $\{0,1\}^n$, and two vertices $x, y \in \{0,1\}^n$ are adjacent in this graph iff they differ in exactly one coordinate.   The next example shows that the hypercube graph belongs to the family of Cayley graphs. 

\begin{Example} \label{eg:hypercube:transpositions} Consider the set of three disjoint transpositions $S = \{ (12), (34), (56)\}$.  Let $H$ be the permutation group generated by $S$.  The Cayley graph $\Cay(H,S)$ is shown in Figure~\ref{fig:cay:3cube}.  This graph is isomorphic to the cube graph $Q_3$.  Observe from the labels of the vertices in the figure that each transposition in $S$ corresponds to one of the dimensions in the edge set of $Q_3$. 
\begin{figure}
\begin{centering}
 \begin{tikzpicture}[scale=0.9,auto]
 \pgfmathsetmacro{\phi}{0.7}
 \pgfmathsetmacro{\delx}{6.5} 
 
  \vertex[fill] (v000) at (0,0) [label=below:{$e$}] {};
  \vertex[fill] (v001) at (2,0) [label=below:{$(12)$}] {};
  \vertex[fill] (v100) at (0,2) [label=left:{$(56)$}] {};
  \vertex[fill] (v101) at (2,2) [label=right:{$(12)(56)$}] {};

  \vertex[fill] (v010) at (0+\phi,0+\phi) {};
  \node at (0.1+\phi+0.35,0+\phi+0.3) {$(34)$};
  \vertex[fill] (v011) at (2+\phi,0+\phi) [label=right:{$(12)(34)$}] {};
  \vertex[fill] (v110) at (0+\phi,2+\phi) [label=above:{$(34)(56)$}] {};
  \vertex[fill] (v111) at (2+\phi,2+\phi) [label=right:{$(12)(34)(56)$}] {};

  \draw (v000) -- (v001) -- (v101) -- (v100) -- (v000);
  \draw (v010) -- (v011) -- (v111) -- (v110) -- (v010);
  \draw (v001) -- (v011);
  \draw (v000) -- (v010);
  \draw (v100) -- (v110);
  \draw (v101) -- (v111);
  \vertex[fill] (w000) at (0+\delx,0) [label=below:{$000$}] {};
  \vertex[fill] (w001) at (2+\delx,0) [label=below:{$100$}] {};
  \vertex[fill] (w100) at (0+\delx,2) [label=left:{$001$}] {};
  \vertex[fill] (w101) at (1.9+\delx,2) [label=right:{$101$}] {};

  \vertex[fill] (w010) at (0+\phi+\delx,0+\phi) {};
  \node at (0+\phi+0.35+\delx,0+\phi+0.3) {$010$};
  \vertex[fill] (w011) at (2+\phi+\delx,0+\phi) [label=right:{$110$}] {};
  \vertex[fill] (w110) at (0+\phi+\delx,2+\phi) [label=above:{$011$}] {};
  \vertex[fill] (w111) at (2+\phi+\delx,2+\phi) [label=right:{$111$}] {};

  \draw (w000) -- (w001) -- (w101) -- (w100) -- (w000);
  \draw (w010) -- (w011) -- (w111) -- (w110) -- (w010);
  \draw (w001) -- (w011);
  \draw (w000) -- (w010);
  \draw (w100) -- (w110);
  \draw (w101) -- (w111);
  
\end{tikzpicture}
\caption{The cube graph as a Cayley graph of a permutation group}
\label{fig:cay:3cube}
\end{centering}
\end{figure}
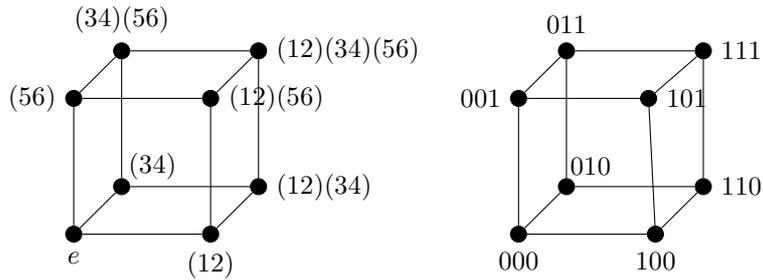
More generally, it can be seen that the Cayley graph of the permutation group generated by the set of $n$ disjoint transpositions $(1,2), (3,4), \ldots, (2n-1,2n)$ is isomorphic to $Q_n$.
\qed
\end{Example}

\begin{Example}
The Cayley graph of $S_3$ with respect to the generating set $S = \{ (12), (23) \}$ is isomorphic to the 6-cycle graph $C_6$ (see Figure~\ref{fig:cay:transp:c6}). The Cayley graph of $S_3$ with respect to the generating set $S = \{ (12), (23), (12) \}$ is the complete bipartite graph $K_{3,3}$ (see Figure~\ref{fig:cay:transp:K33}).
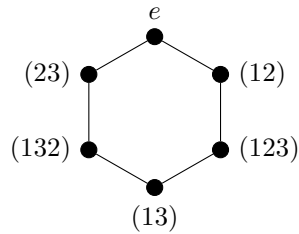
\begin{figure}
\begin{centering}
\begin{tikzpicture}[scale=1,auto]
\vertex[fill] (v1) at (0,1) [label=above:$e$] {};
\vertex[fill] (v2) at (0.866, 0.5) [label=right:$(12)$] {};
\vertex[fill] (v3) at (0.866, -0.5) [label=right:$(123)$] {};
\vertex[fill] (v4) at (0,-1) [label=below:$(13)$] {};
\vertex[fill] (v5) at (-0.866, -0.5) [label=left:$(132)$] {};
\vertex[fill] (v6) at (-0.866,0.5) [label=left:$(23)$] {};
\draw (v1) -- (v2) -- (v3) -- (v4) -- (v5) -- (v6) -- (v1);
\end{tikzpicture}
\caption{The graph $\Cay(S_3, \{ (12), (23) \})$} 
\label{fig:cay:transp:c6}
\end{centering} 
\end{figure}

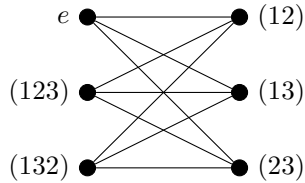
\begin{figure}
\begin{centering}
\begin{tikzpicture}[scale=1,auto]
  \vertex[fill] (x1) at (0,2) [label=left:{$e$}] {};
  \vertex[fill] (x2) at (0,1) [label=left:{$(123)$}] {};
  \vertex[fill] (x3) at (0,0) [label=left:{$(132)$}] {};
  \vertex[fill] (y1) at (2,2) [label=right:{$(12)$}] {};
  \vertex[fill] (y2) at (2,1) [label=right:{$(13)$}] {};
  \vertex[fill] (y3) at (2,0) [label=right:{$(23)$}] {};

\draw (x1) -- (y1);
\draw (x1) -- (y2);
\draw (x1) -- (y3);
\draw (x2) -- (y1);
\draw (x2) -- (y2);
\draw (x2) -- (y3);
\draw (x3) -- (y1);
\draw (x3) -- (y2);
\draw (x3) -- (y3);
\end{tikzpicture}
\caption{The graph $\Cay(S_3, \{ (12), (23) , (13)\})$} 
\label{fig:cay:transp:K33}
\end{centering} 
\end{figure}
If $S$ is a set of transpositions in $S_n$, then the Cayley graph $\Cay(S_n,S)$ is bipartite, with the bipartition consisting of the set of even permutations and the set of odd permutations in $S_n$. 
\qed
\end{Example}

A {\em circulant graph}, denoted by $C_n \langle a_1, \ldots, a_k\rangle$, is defined to be the Cayley graph of the cyclic group $\mathbb{Z}_n = \{0,1,\ldots,n-1\}$ (where the group operation is addition modulo $n$) with respect to the generating set $\{ \pm a_1, \ldots, \pm a_k \}$.  We assume $0 < a_1 < \cdots < a_k < \frac{n+1}{2}$.   This graph can be constructed by placing $n$ points $0,1,\ldots,n-1$ on a circle and joining vertices $i$ and $j$ whenever the absolute difference $|i-j|$ is equal to $a_i$ for some $i$.  If $n$ is even and the last generator $a_k$ is equal to $n/2$, then this generator is called a {\em diagonal jump} since it joins vertices that are diametrically opposite on the circle.   Observe that the circulant $C_n \langle 1 \rangle$ is the Cayley graph $\Cay(\mathbb{Z}_n, \{ \pm 1 \})$, which is the cycle graph on $n$ vertices.  If $X$ is a circulant graph, then its adjacency matrix is a circulant matrix. If $k \ge 2$, then the circulant graph has a 4-cycle $(0, a_1, a_1+a_2, a_2)$.  Since the Petersen graph does not contain 4-cycles, the Petersen graph is not a circulant graph. 

An outstanding open problem in the area of Cayley graphs is the following conjecture:  if $X$ is a finite, connected Cayley graph on 3 or more vertices, then $X$ has a Hamilton cycle.  It is known that if $X$ is a connected, vertex-transitive graph on $n$ vertices, then $X$ has a cycle of length at least $\sqrt{3n}$.   

Standard texts which discuss Cayley graphs include Bollob\'as \cite{Bollobas:1998}, Godsil and Royle \cite{Godsil:Royle:2001} and Biggs \cite{Biggs:1993}. Expository and survey articles on Cayley graphs include  Heydemann \cite{Heydemann:1997}, Alspach \cite{Alspach:1997}, Konstantinova \cite{Konstantinova:2008} and Lakshmivarahan et al \cite{Lakshmivarahan:etal:1993}.  

\subsection{Interconnection networks}
\label{sec:intro:icn}

In this section, we give a brief overview of the area of interconnection networks and their performance measures. The hypercubes, star graphs, degree-diameter problem, cartesian products of graphs and torus graphs are also discussed.  

Modern supercomputers have a large number of processors working in parallel.  One of the fastest computers in the world today is Fujitsu's K-Computer, which has about 500,000 processors.  An interconnection network is a network of interconnected devices and refers to the network used to route data between the processors in a multiprocessor computing system.  The interconnection network is typically modeled as a graph.  The vertices of the graph correspond to  processors, and two vertices are adjacent in the graph whenever there is a direct communication link between the two corresponding processors.  In the sequel, we use the terms interconnection networks and graphs interchangeably.  

One measure of the performance of an interconnection network is its diameter.  The diameter of a graph is defined as the maximum distance between two vertices of the graph.  It's clear the diameter represents the maximum possible delay for communication between two nodes in the interconnection network.  The hypercube graph $Q_n$ has $2^n$ vertices and its diameter is equal to $n$.  Some parallel computers designed in the 1980s used the hypercube topology for the interconnection network.  Modern supercomputers have a very large number of processors (up to a million processors) and use low-dimensional topologies such as the torus graphs. 

The star graph was proposed as a topology of interconnection networks by Akers and Krishnamurthy \cite{Akers:Krishnamurthy:1989}.  The {\em star graph $ST_n$} is defined to be the Cayley graph of the symmetric group $S_n$ with respect to the generating set $S = \{(12), (13), \ldots, (1n)\}$.   Thus, the star graph $ST_n$ is a vertex-transitive graph on $n!$ vertices, with each vertex having degree $n-1$. 

\begin{Theorem}
The diameter of the star graph $ST_n$ ($n \ge 4$) is $\lfloor 3(n-1)/2 \rfloor$.  
\end{Theorem}

Observe that the diameter of the star graph is sublogarithmic in the number of its vertices, making the star graph a superior topology compared to the hypercube from the standpoint of the diameter measure.   

As Table~\ref{table:Qn:STn:degree:diameter} shows, the star graph has better degree-diameter properties than the hypercube in the sense that for a given degree and diameter, the star graph has more vertices. In fact, $ST_8$ has 40 times as many vertices of $Q_{10}$ despite having the same diameter and a smaller degree.  Therefore, a computing system which uses the star graph (rather than hypercube) as its interconnection network can have more processors without sacrificing performance.   The latest results on the degree/diameter problem for Cayley graphs are on the Combinatorics Wiki website \cite{CombWiki}; for maximum degree $\Delta=7$ and diameter $D=10$, there exists a Cayley graph having $6,007,230$ vertices. This graph is much larger than the star graph $ST_8$ (which, for the same degree and diameter, has only $40,320$ vertices).  

The construction of large networks (i.e., networks having a large number of vertices) with given constraints on the maximum degree and diameter is a classical problem in graph theory known as the degree-diameter problem.   A graph with maximum degree $\Delta$ and diameter $D$ is called a $(\Delta,D)$ graph and can have at most $1+ \Delta \sum_{i=0}^{D-1} (\Delta-1)^i$ vertices since the breadth-first search tree starting from a vertex has at most $\Delta$ vertices in the first layer and at most $\Delta(\Delta-1)^{i-1}$ vertices in the $i$th layer.  This upper bound on the maximum possible number of vertices of a $(\Delta,D)$ graph is known as the {\em Moore bound}. An open problem is to construct $(\Delta,D)$ graphs whose number of vertices is as large as possible.  Many of the current record holders happen to be Cayley graphs  \cite{Campbell:etal:1992} \cite{Miller:Siran:2013} \cite{Exoo:Jajcay:2013}.  There is still much room for improvement because the number of vertices in some of these record holders is less than 10\% of the theoretical upper bound (the Moore bound). For example, the maximum possible number of vertices (Moore bound) for a graph with degree 7 and diameter 10 is $1+7 \sum_{i=0}^{i=9}6^i = 84,652,646$, which is more than 10 times the number of vertices in the Cayley graph mentioned above. 
There are also restricted versions of the degree-diameter problem, wherein one may impose additional restrictions on the graph.  For example, one can require the graph to have high symmetry; a general open problem is the following:  given parameters $\Delta$ and $D$, construct large graphs whose maximum degree is at most $\Delta$, diameter is at most $D$, and such that the graph is vertex-transitive.

\begin{table}
\begin{center}
\begin{tabular}{|c|c|c|c|}
\hline
 graph $X$ & $|V|$ & degree $\Delta$ & diameter $D$  \\ \hline
 $Q_n$ & $2^n$ & $n$ & $n$ \\ \hline
 $ST_n$ & $n!$ & $n-1$ & $\lfloor 3(n-1)/2 \rfloor$ \\ \hline
 $Q_{10}$ & $1,024$ & 10 & 10 \\  \hline
 $ST_8$ & $40,320$ & 7 & 10 \\ \hline
  \end{tabular}
\caption{Degree-diameter properties of the hypercube $Q_n$ and star graph $ST_n$.}
\end{center}
 \end{table}
\label{table:Qn:STn:degree:diameter}

The fastest supercomputers (today) use the torus interconnection network, which can be defined as the cartesian product of cycles (or as a Cayley graph, as we shall see).  

\begin{Definition}
 Let $Y=(V,E)$ and $Y'=(V',E')$ be two graphs.  The {\em cartesian product} of $Y$ and $Y'$, denoted by $Y \square Y'$,  is the graph with vertex set $V \times V'$, and two vertices $(u,u'), (v, v') \in V \times V'$ are adjacent in $Y \square Y'$ whenever $u=v$ and $u' v' \in E(Y')$ or $u'=v'$ and $uv \in E(Y)$.  
\end{Definition}

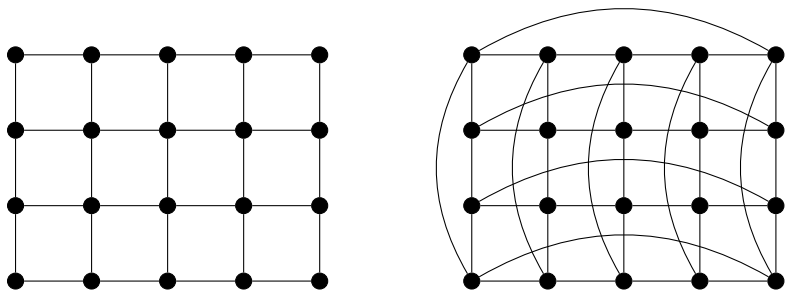
\begin{figure}
\begin{centering}
\begin{tikzpicture}[scale=1,auto]

\vertex[fill] (v00) at (0,0) {};
\vertex[fill] (v01) at (0,1) {};
\vertex[fill] (v02) at (0,2) {};
\vertex[fill] (v03) at (0,3) {};
\draw (v00) -- (v01) -- (v02) -- (v03);

\vertex[fill] (v10) at (1,0) {};
\vertex[fill] (v11) at (1,1) {};
\vertex[fill] (v12) at (1,2) {};
\vertex[fill] (v13) at (1,3) {};
\draw (v10) -- (v11) -- (v12) -- (v13);

\vertex[fill] (v20) at (2,0) {};
\vertex[fill] (v21) at (2,1) {};
\vertex[fill] (v22) at (2,2) {};
\vertex[fill] (v23) at (2,3) {};
\draw (v20) -- (v21) -- (v22) -- (v23);

\vertex[fill] (v30) at (3,0) {};
\vertex[fill] (v31) at (3,1) {};
\vertex[fill] (v32) at (3,2) {};
\vertex[fill] (v33) at (3,3) {};
\draw (v30) -- (v31) -- (v32) -- (v33);

\vertex[fill] (v40) at (4,0) {};
\vertex[fill] (v41) at (4,1) {};
\vertex[fill] (v42) at (4,2) {};
\vertex[fill] (v43) at (4,3) {};
\draw (v40) -- (v41) -- (v42) -- (v43);

\draw (v00) -- (v10) -- (v20) -- (v30) -- (v40);
\draw (v01) -- (v11) -- (v21) -- (v31) -- (v41);
\draw (v02) -- (v12) -- (v22) -- (v32) -- (v42);
\draw (v03) -- (v13) -- (v23) -- (v33) -- (v43);


\vertex[fill] (w00) at (6,0) {};
\vertex[fill] (w01) at (6,1) {};
\vertex[fill] (w02) at (6,2) {};
\vertex[fill] (w03) at (6,3) {};
\draw (w00) -- (w01) -- (w02) -- (w03);

\vertex[fill] (w10) at (7,0) {};
\vertex[fill] (w11) at (7,1) {};
\vertex[fill] (w12) at (7,2) {};
\vertex[fill] (w13) at (7,3) {};
\draw (w10) -- (w11) -- (w12) -- (w13);

\vertex[fill] (w20) at (8,0) {};
\vertex[fill] (w21) at (8,1) {};
\vertex[fill] (w22) at (8,2) {};
\vertex[fill] (w23) at (8,3) {};
\draw (w20) -- (w21) -- (w22) -- (w23);

\vertex[fill] (w30) at (9,0) {};
\vertex[fill] (w31) at (9,1) {};
\vertex[fill] (w32) at (9,2) {};
\vertex[fill] (w33) at (9,3) {};
\draw (w30) -- (w31) -- (w32) -- (w33);

\vertex[fill] (w40) at (10,0) {};
\vertex[fill] (w41) at (10,1) {};
\vertex[fill] (w42) at (10,2) {};
\vertex[fill] (w43) at (10,3) {};
\draw (w40) -- (w41) -- (w42) -- (w43);

\draw (w00) -- (w10) -- (w20) -- (w30) -- (w40);
\draw (w01) -- (w11) -- (w21) -- (w31) -- (w41);
\draw (w02) -- (w12) -- (w22) -- (w32) -- (w42);
\draw (w03) -- (w13) -- (w23) -- (w33) -- (w43);

\path[bend left] (w00) edge (w40);
\path[bend left] (w01) edge (w41);
\path[bend left] (w02) edge (w42);
\path[bend left] (w03) edge (w43);

\path[bend left] (w00) edge (w03);
\path[bend left] (w10) edge (w13);
\path[bend left] (w20) edge (w23);
\path[bend left] (w30) edge (w33);
\path[bend left] (w40) edge (w43);

\end{tikzpicture}
\caption{The $4 \times 5$ mesh and the $4 \times 5$ torus.}
\label{fig:torus:4x5}
\end{centering}
\end{figure}

Let $P_n$ and $C_n$ denote the path graph and cycle graph, respectively, on $n$ vertices. The cartesian product $P_4 \square P_5$ is the {\em $4 \times 5$ mesh graph} shown in the first part of Figure~\ref{fig:torus:4x5}.  The cartesian product $C_4 \square C_5$ is the {\em $4 \times 5$ torus graph} shown in the second part of Figure~\ref{fig:torus:4x5}.  Observe from the figure that the $4 \times 5$ torus graph can be obtained from the $4 \times 5$ mesh graph by adding so-called wrap-around edges between the first node and last node in each row (and likewise in each column). These edges reduce the diameter to about half its value since the diameter of of the $n$-cycle graph is about half the diameter of the path graph on $n$ vertices.   The definition of the cartesian product of two graphs can be extended to the cartesian product of any number of graphs \cite{Imrich:Klavzar:Rall:2008}.  The hypercube $Q_n$ is isomorphic to the cartesian product of $n$ copies of the complete graph $K_2$. 

The cycle graph $C_r$ is isomorphic to the Cayley graph $\Cay(\mathbb{Z}_r, \{ \pm 1 \})$, where $\mathbb{Z}_r$ is the cyclic group consisting of the elements $\{0,1,\ldots,r-1\}$ and with the group operation in $\mathbb{Z}_r$ being addition modulo $r$.  Observe that the 2D-torus graph $C_4 \square C_5$ shown in Figure~\ref{fig:torus:4x5} is isomorphic to the Cayley graph $\Cay(\mathbb{Z}_4 \times \mathbb{Z}_5, S)$, where the generating set is $S = \{\pm (1,0), \pm(0,1) \}$.  More generally,

\begin{Definition}
A {\em $k$-dimensional torus graph} is defined to be a graph which is the cartesian product of $k$ cycle graphs.   
\end{Definition}

Hence, a $k$-dimensional torus graph is isomorphic to the Cayley graph of some group $\mathbb{Z}_{r_1} \times \cdots \times \mathbb{Z}_{r_k}$ with respect to the generating set $\{\pm e_1, \ldots, \pm e_k\}$, where $e_i$ consists of a 1 in the $i$th coordinate and 0 in the remaining coordinates.   The torus is said to be {\em balanced} if all the $r_i$'s are equal.  Torus graphs which are Cayley graphs of the group $\mathbb{Z}_r^k$ are also called $r$-ary $k$-cubes.

Since the topology of a $k$-dimensional torus interconnect is a cartesian product, each node is a $k$-tuple and routing is simple. For example, in the $6 \times 5 \times 6$ torus graph, the vertices are the elements of $\mathbb{Z}_6 \times \mathbb{Z}_5 \times \mathbb{Z}_6$,  and a shortest path from node $(1,0,3)$ to node $(3, 4, 1)$ is easily determined to be via nodes $(2,0,3)$, $(3,0,3)$, $(3,4,3)$, $(3,4,2)$, $(3,4,1)$.   The diameter of this torus is $3 + 2 + 3 = 8$.  

Suppose we need a topology for a computing system with about 1000 nodes.  A $32 \times 32$ torus (which has 2 dimensions) has diameter $16+16=32$.  A $10 \times 10 \times 10$ torus has diameter $5+5+5=15$.  It is clear that a balanced torus is better than an unbalanced one (such as a $2 \times 10 \times 50$ torus).  If the number of dimensions increases, the diameter reduces further - for example, the hypercube $Q_{10}$ is a 10-dimensional torus and has diameter 10.  

The diameter is only one measure of the performance of the interconnection network.  A realistic comparison of topologies must take into account the wire lengths and various other parameters. Dally \cite{Dally:1990} showed, using models for delay, that low-dimensional networks (such as tori) have lower latency than high-dimensional networks (such as $n$-dimensional hypercubes).  The torus interconnect has the advantage that its logical graph and physical graph are identical - the physical network of a torus interconnect is such that two nodes which are adjacent in the torus graph are also physically close to each other. While the hypercube has a smaller diameter than the torus, in practice the hypercube is more expensive because it is more difficult to physically lay out or build a high-dimensional hypercube interconnect in the physical form of a 2D or 3D array.  The fastest supercomputers today use a torus interconnect - IBM's Blue Gene/Q uses a 5D torus, and Fujitsu's K Computer uses a 6D torus \cite{Ajima:etal:2012}.  

\section{Automorphism groups of Cayley graphs}
\label{sec:aut:cayley}

In this section, we recall some fundamental results on the automorphism groups of Cayley graphs.  We show that the automorphism group of every Cayley graph $\Cay(H,S)$ contains the following two subgroups: the right regular representation $R(H)$ and the set $\Aut(H,S)$ of automorphisms of $H$ that fixes $S$ setwise.  Then, we define normal Cayley graphs, obtain the automorphism group of the hypercube, and survey some recent results in the literature.

Given a group $H$ and a subset $S \subseteq H$, the problem of obtaining the automorphism group of the Cayley graph $\Cay(H,S)$ is in general a difficult, open and interesting problem \cite{Godsil:1981} 
\cite{Xu:1998}
\cite{Jajcay:1994}
\cite{Jajcay:2000}.  The automorphism groups of various families of Cayley graphs have been obtained recently.  The automorphism groups of Cayley graphs of the symmetric group (for various generating sets) are investigated in 
\cite{Zhang:Huang:2005}
\cite{Feng:2006} 
\cite{Steinhardt:2007}
\cite{Ganesan:DM:2013} 
\cite{Ganesan:JACO} 
\cite{Ganesan:AJC:2016}  
\cite{Ganesan:DMGT:toappear}
\cite{Deng:Zhang:2011}
\cite{Deng:Zhang:2012}.
The automorphism group of the hypercube and its supergraphs (which are Cayley graphs of $\mathbb{Z}_2^n$) are studied in 
\cite{Harary:2000}
\cite{Choudum:Sunitha:2008}
\cite{Mirafzal:2011}.
Other papers on this topic include  
\cite{Zhou:2011} 
\cite{Dobson:2012}.

\subsection{All Cayley graphs are vertex-transitive}

Let $H$ be a group. Throughout this section, we assume that $S$ is a subset of a finite group $H$ satisfying the condition $e \notin S = S^{-1}$.  Thus, the Cayley graph $\Cay(H,S)$ is a simple, undirected graph.  

Given a group $H$, the right regular representation of $H$, denoted by $R(H)$, is the permutation group $\{r_h: h \in H\}$ in $\Sym(H)$, where $r_h$ is defined to be the map $x \mapsto xh$ from $H$ to itself.  Recall that Cayley's theorem asserts that every group $H$ is isomorphic to the permutation group $R(H)$. 

\begin{Theorem}
Let $H$ be a group and let $S$ be a subset of $H$. Let $X$ be the Cayley graph $\Cay(H,S)$.  Then, the right regular representation $R(H)$ is a subgroup of $\Aut(X)$.   
\end{Theorem}

\noindent \emph{Proof:}
Let $u, v \in H$ be vertices of the Cayley graph $X=\Cay(H,S)$. Let $r_h \in R(H)$.  Then, $v$ and $u$ are adjacent in $X$ if and only if $v = su$ for some $s \in S$, iff $vh = suh$ for some $s \in S$, iff $vh$ and $uh$ are adjacent in $X$.  Thus, $r_h$ is an automorphism of $X$. 
\qed

Since $R(H)$ acts transitively on $H$, we get the following result:

\begin{Corollary}
All Cayley graphs are vertex-transitive.
\end{Corollary}

Since $R(H)$ acts regularly on $H$, the automorphism group of every Cayley graph $\Cay(H,S)$ contains a regular subgroup.   Sabidussi \cite{Sabidussi:1958} proved the following converse. 

\begin{Theorem}
If the automorphism group of a graph $X$ contains a regular subgroup, then $X$ is a Cayley graph.
\end{Theorem}

\noindent \emph{Proof:}
Let $G$ be a regular subgroup in $\Aut(X)$.  We show that there exists a group $H$ and a subset $S \subseteq H$ such that $X \cong \Cay(H,S)$.  Let $V(X) = \{v_1,\ldots,v_n\}$. Since $G$ acts regularly on $V(X)$, $G$ is of the form $G = \{g_1,\ldots,g_n\}$, where $g_i$ is the unique element of $G$ which takes $v_1$ to $v_i$.  In particular, $g_1$ is the identity element in $G$.  Let $H := \{g_1,\ldots,g_n\}$, and let 
\[S := \{ g_i \in G: v_i \mbox{ is adjacent to } v_1 \mbox{ in } X \}. \]
We claim that $e \notin S, S = S^{-1}$ and $X \cong \Cay(H,S)$.  Since $X$ has no self-loops, $v_1$ is not adjacent to itself in $X$. Hence, $g_1 \notin S$.  Suppose $g_i \in S$. Then $\{v_1,v_i\} = \{v_1, v_1^{g_i} \}\in E(X)$.  Since $g_i^{-1}$ is an automorphism of $X$, $\{v_1^{g_i^{-1}}, v_1\} \in E(X)$.  This implies that $g_i^{-1} \in S$, whence $S$ is closed under inverses.  Finally, we show that $v_i$ and $v_j$ are adjacent in $X$ iff $g_i$ and $g_j$ are adjacent in $\Cay(H,S)$.  We have the following sequence of statements: $v_i$ and $v_j$ are adjacent in $X$ iff $v_1^{g_i}$ and $v_1^{g_j}$ are adjacent in $X$, iff $v_1$ and $v_1^{g_j g_i^{-1}}$ are adjacent in $X$, iff $g_j g_i^{-1} \in S$, iff $g_i$ and $g_j$ are adjacent in $\Cay(H,S)$. 
\qed

While every Cayley graph is vertex-transitive, the converse is not true in general  - there exist vertex-transitive graphs which are not Cayley graphs: 

\begin{Theorem} \label{thm:Petersen:not:cayley}
 The Petersen graph is not a Cayley graph.
\end{Theorem}

\noindent \emph{Proof}: The Petersen graph is a 3-regular, connected graph on 10 vertices.  There are exactly two groups of order 10, namely the cyclic $C_{10}$ and the dihedral group $D_{10}$.  If the Petersen graph $X$ is a Cayley graph, then it is isomorphic to $\Cay(C_{10},S)$ or $\Cay(D_{10},S)$ for some $S$ such that $S=S^{-1}$, $|S|=3$ and $\langle S \rangle=C_{10}$ or $\langle S \rangle = D_{10}$. Since the girth of the Petersen graph is 5, it suffices to show that every cubic, connected Cayley graph on 10 vertices has a 4-cycle.  

\begin{figure}
\begin{centering}
\begin{tikzpicture}[scale=1,auto]
\vertex[fill] (u1) at (2,2) [label=right:$e$] {};
\vertex[fill] (u2) at (0,2) [label=left:$a$] {};
\vertex[fill] (u3) at (0,0) [label=left:$ba$] {};
\vertex[fill] (u4) at (2,0) [label=right:$b$] {};
\draw [->] (u1) edge node {$a$} (u2);
\draw [->] (u2) edge node {$b$} (u3);
\draw [->] (u3) edge node {$a^{-1}$} (u4);
\draw [->] (u4) edge node {$b^{-1}$} (u1);

\vertex[fill] (v1) at (6,1) [label=left:$e$] {};
\vertex[fill] (v2) at (8,0) [label=below:$b$] {};
\vertex[fill] (v3) at (8,2) [label=above:$a$] {};
\vertex[fill] (v4) at (10,1) [label=right:{$ab=ba$}] {};
\draw [-] (v1) edge (v2);
\draw [-] (v1) edge (v3);
\draw [-] (v2) edge node {$a$} (v4);
\draw [-] (v3) edge node {$b$} (v4);
\end{tikzpicture}
\caption{Cayley graphs of abelian groups have 4-cycles.}
\label{fig:abelian:cay:have:4cycles}
\end{centering}
\end{figure}
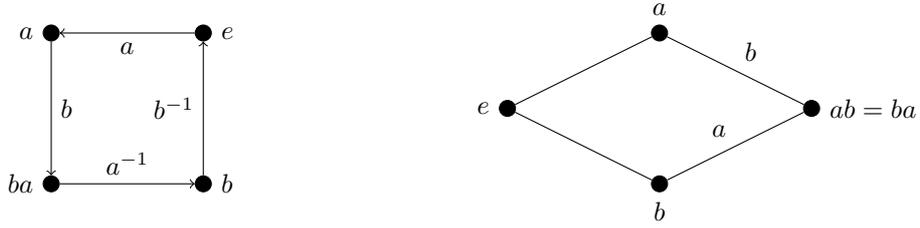

Let $H$ be an abelian group, and suppose $1 \notin S=S^{-1} \supseteq \{a,a^{-1}, b \}$.  We claim that $\Cay(H,S)$ contains a 4-cycle.  Since $H$ is abelian, $a^{-1}ba=b$, and so the sequence of vertices $e,a,ba,a^{-1}ba=b,e$ forms a 4-cycle in $\Cay(H,S)$ (see Figure~\ref{fig:abelian:cay:have:4cycles}). This proves the claim.  Taking $H=C_{10}$, we get that the Petersen graph is not a Cayley graph of the cyclic group $C_{10}$.

The dihedral group $D_{10}$ has the presentation $\langle r,s | r^5=s^2=e, rs=sr^{-1} \rangle$.  Suppose $S \subseteq D_{10}$, $S=S^{-1}$ and $|S|=3$.  We show that $\Cay(D_{10},S)$ has a 4-cycle.   We consider two cases for the 3-element generator set $S=S^{-1}$: either $S$ is of the form $\{a,a^{-1},b\}$ or $S$ is of the form $\{a,b,c\}$ consisting of three distinct elements of order 2.  

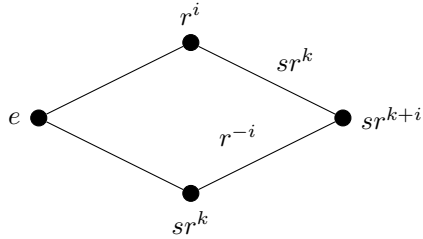
\begin{figure}
\begin{centering}
\begin{tikzpicture}[scale=1,auto]
\vertex[fill] (v1) at (6,1) [label=left:$e$] {};
\vertex[fill] (v2) at (8,0) [label=below:$sr^k$] {};
\vertex[fill] (v3) at (8,2) [label=above:$r^i$] {};
\vertex[fill] (v4) at (10,1) [label=right:{$sr^{k+i}$}] {};
\draw [-] (v1) edge (v2);
\draw [-] (v1) edge (v3);
\draw [-] (v2) edge node {$r^{-i}$} (v4);
\draw [-] (v3) edge node {$sr^k$} (v4);
\end{tikzpicture}
\caption{The Cayley graph $\Cay(D_{10},S)$ contains a 4-cycle when  $S=\{r^i, r^{-i}, sr^k \}$.}
\label{fig1:cay:D10:has:4cycles}
\end{centering}
\end{figure}

\begin{figure}
\begin{centering}
\begin{tikzpicture}[scale=0.5,auto]
\vertex[fill] (v1) at (6,1) [label=left:$e$] {};
\vertex[fill] (v2) at (8,0) [label=below:$sr^j$] {};
\vertex[fill] (v3) at (8,2) [label=above:$sr^i$] {};
\vertex[fill] (v3b) at (8,-2) [label=below:$sr^k$] {};
\vertex[fill] (v4) at (10,1) [label=right:{$r^{j-i}$}] {};
\draw [-] (v1) edge (v2);
\draw [-] (v1) edge (v3);
\draw [-] (v1) edge (v3b);
\draw [-] (v2) edge (v4);
\end{tikzpicture}
\caption{The Cayley graph $\Cay(D_{10},S)$ contains a 4-cycle when $S=\{sr^i, sr^j, sr^k\}$.}
\label{fig2:cay:D10:has:4cycles}
\end{centering}
\end{figure}
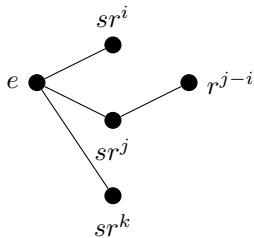

First suppose $S$ is of the form $S=\{r^i, r^{-i}, sr^k \}$.  Then $r^i$ and $sr^k$ have $sr^{k+i}$ as a common neighbor because $sr^k r^i=r^{-i}sr^k$. This gives rise to a 4-cycle (see Figure~\ref{fig1:cay:D10:has:4cycles}) in the Cayley graph.  Now suppose $S$ is of the form $S = \{sr^i, sr^j, sr^k\}$, with $i,j,k \in \{0,1,2,3,4\}$ distinct. We have $sr^i sr^j = r^{j-i}$ (see Figure~\ref{fig2:cay:D10:has:4cycles}). Consider the six differences $i-j, j-i, i-k, k-i, j-k$ and $k-j$ that take values in $\{0,1,2,3,4\}$.  
By the pigeonhole principle, some two differences are equal, giving rise to a 4-cycle in the Cayley graph $\Cay(D_{10},S)$. 
\qed

In the proof of Theorem~\ref{thm:Petersen:not:cayley}, we also showed the following:

\begin{Corollary}
 Let $H$ be an abelian group and suppose $1 \notin S=S^{-1} \subseteq H$ and $|S| \ge 3$.  Then, the girth of the Cayley graph $\Cay(H,S)$ is at most 4.
\end{Corollary}

We saw that not every connected vertex-transitive graph is a Cayley graph. It turns out that every connected vertex-transitive graph is a retract of some Cayley graph \cite{Sabidussi:1964}.  Another way to state this result is: Cayley graphs represent all vertex-transitive graphs up to homomorphism-equivalence.  Equivalently, every vertex-transitive graph is isomorphic to some Cayley coset graph.

The automorphism group $G:=\Aut(\Cay(H,S))$ of the Cayley graph $\Cay(H,S)$ can be expressed as $G = R(H) G_e := \{r_h g: r_h \in R(H), g \in G_e\}$, where $G_e$ is the stabilizer in $G$ of the identity vertex $e$.  The automorphism group $G$ of every Cayley graph $\Cay(H,S)$ can be written as a so-called rotary product $G = R(H) \times_{rot} G_e$ (cf. \cite{Jajcay:1994}, \cite{Jajcay:2000}), which is a generalization of a semidirect product.  

Cayley graphs $\Cay(H,S)$ which have the smallest possible full automorphism group $R(H)$ have been investigated further:.

\begin{Definition}
 A Cayley graph $\Cay(H,S)$ is said to be a {\em graphical regular representation} (or {\em GRR}) if its full automorphism group is equal to $R(H)$. 
\end{Definition}

An example of an infinite family of GRR's is the family of Cayley graphs of the symmetric group $S_n$ with respect to transposition sets $S$ such that the transposition graph is an asymmetric tree (cf. Theorem~\ref{thm:GR:autgroup:asym:tree}).  An open conjecture in the literature is the following: if $H$ does not belong to an exceptional family of groups, then almost all Cayley graphs of $H$ are GRRs \cite{Babai:Godsil:1982}.   

\subsection{Normal Cayley graphs}

Normal Cayley graphs are Cayley graphs that have the smallest possible full automorphism group in a certain sense.  In this section, we define normal Cayley graphs.  We begin by recalling some properties of the automorphism group of a group $H$, denoted by $\Aut(H)$. 

\begin{Lemma} \label{lem:AutH:normalizes:RH}
 Let $H$ be a group.  Then, $\Aut(H)$ normalizes $R(H)$. 
\end{Lemma}

\noindent \emph{Proof:}
Let $\sigma \in \Aut(H)$ and let $(r_z: h \mapsto hz) \in R(H)$. We show that $\sigma^{-1} r_z \sigma \in R(H)$.  Let $y$ be any element in $H$ and let $x:=y^{\sigma^{-1}}$.  Then, $\sigma^{-1} r_z \sigma$ takes $y$ to $y^{\sigma^{-1} r_z \sigma} = x^{r_z \sigma} = (xz)^\sigma$ $= x^\sigma z^\sigma =  yz^{\sigma} = y^{r_{(z^\sigma)}}$. Thus, $\sigma^{-1} r_z \sigma = r_{z^\sigma} \in R(H)$. 
\qed

\begin{Definition}
  Let $H$ be a group. 
 The {\em holomorph} of $H$, denoted by $\Hol(H)$, is the normalizer in $\Sym(H)$ of the right regular representation $R(H)$:
 \[ \Hol(H) := N_{\Sym(H)} (R(H)). \]
\end{Definition}

It can be shown that the holomorph of $H$ contains copies of both $H$ and its automorphism group $\Aut(H)$:

\begin{Lemma} \label{lem:holomorph:eq:RHAutH}
 Let $H$ be a group. Then 
 \[N_{\Sym(H)} (R(H)) = R(H) \Aut(H).\]
\end{Lemma}

Thus, $\Hol(H) = R(H) \Aut(H)$, which by  Lemma~\ref{lem:AutH:normalizes:RH} equals $R(H) \rtimes \Aut(H)$.   

Given a group $H$ and a subset $S \subseteq H$, let $\Aut(H,S)$ denote the set of automorphisms of the group $H$ that fixes $S$ setwise: \[\Aut(H,S) := \{g \in \Aut(H): S^g=S\}.\]  We now show that every Cayley graph admits a certain subgroup of automorphisms: 

\begin{Theorem} \label{thm:AutHS:subgroup:of:Ge}
Let $H$ be a group and let $S$ be a subset of $H$.   Then, the set $\Aut(H,S)$ of automorphisms of $H$ that fixes $S$ setwise is a subgroup of the stabilizer $\Aut(\Cay(H,S))_e$. 
\end{Theorem}

\noindent \emph{Proof:} Let $X$ be the Cayley graph $\Cay(H,S)$ and let $G :=\Aut(X)$.  We need to show that $\Aut(H,S) \le G_e$. Let $g \in \Aut(H,S)$.  Since $g$ is an automorphism of the group $H$, $g$ fixes the identity element $e$.  To show that $g$ is an automorphism of $X$, observe that vertices $u$ and $v$ are adjacent in $X$ iff $uv^{-1} \in S$, iff $(uv^{-1})^g \in S$, iff $u^g$ and $v^g$ are adjacent in $X$.  \qed

We showed above that every Cayley graph possesses certain symmetries: the automorphism group of every Cayley graph $\Cay(H,S)$ contains the subgroups $R(H)$ and $\Aut(H,S)$.  Note that $R(H)$ acts regularly on $H$, and every element in $\Aut(H,S)$ fixes the identity element $e$. Thus, $R(H) \cap \Aut(H,S) = 1$.  By Lemma~\ref{lem:AutH:normalizes:RH}, $R(H) \Aut(H,S)$ can be expressed as the semidirect product $R(H) \rtimes \Aut(H,S)$.   Thus, we have the following:

\begin{Corollary} \label{cor:autcay:contains:sdp}
 Let $H$ be a group and let $S$ be a subset of $H$.  Then, the automorphism group of the Cayley graph $\Cay(H,S)$ contains  $R(H) \rtimes \Aut(H,S)$. 
\end{Corollary}

If $\Cay(H,S)$ has the smallest possible full automorphism group in the sense that it contains no other automorphisms besides those in $R(H) \rtimes \Aut(H,S)$, then the Cayley graph $\Cay(H,S)$ is said to be normal:

\begin{Definition}
 A Cayley graph $\Cay(H,S)$ is said to be {\em normal} if its full automorphism group is equal to $R(H) \rtimes \Aut(H,S)$. 
\end{Definition}

We now prove two preliminary lemmas before proving some equivalent conditions for normality.

\begin{Lemma} \label{lem:AutHcapGe:subset:AutHS}
 Let $G$ be the automorphism group of the Cayley graph $\Cay(H,S)$.  Then, $\Aut(H) \cap G_e = \Aut(H,S)$. 
\end{Lemma}

\noindent \emph{Proof:}
If $g \in G_e$, then $g$ is an automorphism of the graph $\Cay(H,S)$ that fixes the identity vertex $e$, whence $g$ also fixes the set $S$ of neighbors of $e$ setwise.  Hence, $\Aut(H) \cap G_e \subseteq \Aut(H,S)$.  The reverse inclusion follows from Theorem~\ref{thm:AutHS:subgroup:of:Ge}.
\qed

\begin{Lemma} \label{lem:NGRH:eq:RHAutHS}
 Let $G$ be the automorphism group of the Cayley graph $\Cay(H,S)$.  Then 
 \[ N_G(R(H)) = R(H) \Aut(H,S). \]
\end{Lemma}

\noindent \emph{Proof:}
We have that $N_G(R(H)) = N_{\Sym(H)} (R(H)) \cap G$, which is equal to $R(H) \Aut(H) \cap G$ by Lemma~\ref{lem:holomorph:eq:RHAutH}.  Hence, $N_G(R(H)) = R(H) \Aut(H) \cap R(H) G_e = R(H) (\Aut(H) \cap G_e) = R(H) \Aut(H,S)$, where the last equality is by Lemma~\ref{lem:AutHcapGe:subset:AutHS}.
\qed

\begin{Theorem} \label{thm:equiv:cond:normality}
 Let $H$ be a group and let $S$ be a subset of $H$.  Let $G$ be the automorphism group of the Cayley graph $\Cay(H,S)$.  The following conditions are equivalent:
 \\ (a) $\Cay(H,S)$ is a normal Cayley graph, i.e. $G=R(H) \rtimes \Aut(H,S)$. 
 \\ (b) $R(H)$ is a normal subgroup of $G$.
 \\ (c) $G_e = \Aut(H,S)$. 
\end{Theorem}

\noindent \emph{Proof:} 
(a) $\implies$ (b): Suppose $G=R(H) \Aut(H,S)$. Then $N_{R(H) \Aut(H,S)} (R(H)) = R(H) \Aut(H,S)$ by Lemma~\ref{lem:NGRH:eq:RHAutHS}.  

(b) $\implies$ (a):  Suppose $R(H)$ is normal in $G$. Then $N_G(R(H)) = G$, and by Lemma~\ref{lem:NGRH:eq:RHAutHS}, $G = R(H) \Aut(H,S)$.  

(a) $\iff$ (c): Since $G = R(H) G_e$ always, $G = R(H) \Aut(H,S)$ iff $\Aut(H,S) = G_e$.  
\qed

Let $G$ be the automorphism group of the Cayley graph $\Cay(H,S)$.  By Theorem~\ref{thm:equiv:cond:normality}, an equivalent condition for the Cayley graph $\Cay(H,S)$ to be non-normal is that $G_e$ be strictly larger than $\Aut(H,S)$.  For in this case, the full automorphism group $G = R(H)  G_e$ is strictly larger than the smallest possible full automorphism group $R(H) \rtimes \Aut(H,S)$. Recall that if $H$ is non-abelian, then the inverse map $(\sigma: H \rightarrow H, h \mapsto h^{-1})$ is not a homomorphism from $H$ to itself and so $\sigma \notin \Aut(H,S)$.  Thus, if the inverse map is an automorphism of the Cayley graph, then the Cayley graph is non-normal.  A family of non-normal Cayley graphs which admit the inverse map as an automorphism is the complete transposition graphs (cf. Theorem~\ref{thm:completeTS:nonnormal} and Theorem~\ref{thm:completeTS:autgroup}).

The problem of obtaining the (full) automorphism group of Cayley graphs is in general an open problem, and determining (for a given $H, S$) whether $\Cay(H,S)$ is a GRR or a normal Cayley graph is an open problem (cf. \cite{Xu:1998}).  

\subsection{Automorphism group of the hypercube}

In this section we obtain the automorphism group of one of the simplest families of Cayley graphs - the hypercube graphs.

Given a graph $X$, let $G=\Aut(X)$. Fix $v \in V$.  Then, $G_v$ denotes the set of automorphisms of $X$ that fixes the vertex $v$, and $L_v = L_v(X)$ denotes the set of automorphisms of $X$ that fixes the vertex $v$ and each of its neighbors.  Thus, $L_v$ is the kernel of the action of $G_v$ on the set $X_1(v)$ of neighbors of $v$.  It follows that $G_v/L_v$ is isomorphic to a subgroup of $\Sym(X_1(v))$.  If the graph $X$ is $k$-regular, then $|G_v| \le |L_v|~(k!)$.  This gives an upper bound on the number of elements in the stabilizer $G_v$. This also gives an upper bound on $|G|$ if the graph is vertex-transitive.  This method is often useful in proving that one has obtained all the automorphisms of a graph.

Recall that the hypercube $Q_n$ is defined to be the graph with vertex set $\{0,1\}^n$, and two vertices are adjacent in $Q_n$ whenever the corresponding vectors differ in exactly one coordinate.  Each $g \in S_n$ acts in a natural way on $V(Q_n)$ by permuting the coordinates of each vector.  If $x$ and $y$ are adjacent vertices in $Q_n$ (i.e. they differ in exactly one coordinate), then $x^g$ and $y^g$ also differ in exactly one coordinate, and conversely.  Thus, the action of $g$ on $V(Q_n)$ preserves adjacency.  It follows that $S_n$ acts on $V(Q_n)$ as a group of automorphisms of $Q_n$.  Also note that if $z$ is any vector, then vertices $x$ and $y$ are adjacent in $Q_n$ iff vertices $x+z$ and $y+z$ are adjacent in $Q_n$; here, addition is performed mod 2 and componentwise. Thus, translation by $z$ is an automorphism of $Q_n$.  We show below that $Q_n$ has no other automorphisms besides these, i.e. that each automorphism of $Q_n$ is the composition of a translation $z \in \{0,1\}^n$ and a permutation $g \in S_n$.  

\begin{Theorem} \label{thm:Qn:autgroup} 
The automorphism group of the hypercube graph $Q_n$ is isomorphic to the semidirect product $\mathbb{Z}_2^n \rtimes S_n$.
\end{Theorem}

\noindent \emph{Proof:}
The finite vector space $\mathbb{F}_2^n$ forms a group $\mathbb{Z}_2^n$ under the operation of vector addition.  The hypercube graph $Q_n$ is isomorphic to the Cayley graph $\Cay(\mathbb{Z}_2^n,S)$, where $S=\{e_1,\ldots,e_n\}$ and $e_i$ is the unit vector consisting of a 1 in the $i$th coordinate and 0 in the remaining coordinates.   By Corollary~\ref{cor:autcay:contains:sdp},  $\Aut(Q_n)$ contains $R(\mathbb{Z}_2^n) \rtimes \Aut(\mathbb{Z}_2^n,S)$.  The set of automorphisms of the group $\mathbb{Z}_2^n$ that fixes $S$ setwise is precisely the set of permutation matrices in the general linear group $GL_n(\mathbb{F}_2)$.  Hence, $\Aut(\mathbb{Z}_2^n,S) \cong S_n$. 

Let $X=Q_n$ and let $G = \Aut(X)$.  We have shown above that $G$ contains a subgroup isomorphic to $\mathbb{Z}_2^n \rtimes S_n$.   It remains to show that the graph $X$ has no other automorphisms besides these $(2^n)n!$ automorphisms. Let $G_e$ be the stabilizer in $G$ of the identity vertex $e \in V(X)$. Let $L_e$ be the set of automorphisms of $X$ that fixes the identity vertex $e$ and each of its neighbors in $X$.  Note that $L_e$ is the kernel of the action of $G_e$ on $S$, whence $|G_e| \le |L_e|~|\Sym(S)|$. 

\begin{figure}
\begin{centering}
\begin{tikzpicture}[scale=1.5,auto]
  \vertex[fill] (v000) at (0,1) [label=left:$e$] {};

  \vertex[fill] (v100) at (2,2) [label=below:$100$] {};
  \vertex[fill] (v010) at (2,1) [label=below:$010$] {};
  \vertex[fill] (v001) at (2,0) [label=below:$001$] {};

  \vertex[fill] (v110) at (4,2) [label=below:$110$] {};
  \vertex[fill] (v101) at (4,1) [label=below:$101$] {};
  \vertex[fill] (v011) at (4,0) [label=below:$011$] {};

  \vertex[fill] (v111) at (6,1) [label=right:$111$] {};

  \draw (v000) -- (v100);
  \draw (v000) -- (v010);
  \draw (v000) -- (v001);

  \draw (v001) -- (v101);
  \draw (v001) -- (v011);
  \draw (v010) -- (v011);
  \draw (v010) -- (v110);
  \draw (v100) -- (v101);
  \draw (v100) -- (v110);

  \draw (v011) -- (v111);
  \draw (v101) -- (v111);
  \draw (v110) -- (v111);

  \node at (2, -1) {$X_1(e)$};
  \node at (4, -1) {$X_2(e)$};
  \node at (6, -0.5) {$X_3(e)$};
\end{tikzpicture}
\caption{The distance partition of $X=Q_3$ with respect to the identity vertex $e$.}
\label{fig:dist:partition:Q3}
\end{centering}
\end{figure}
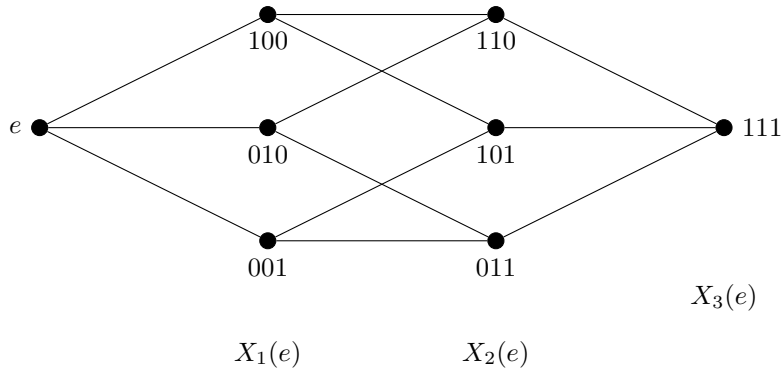

We claim that $L_e=1$.  Let $g \in L_e$.  Then, $g$ fixes the identity vertex $e$ and each of its neighbors.  Let $X_i(e)$ denote the set of vertices of $X$ whose distance to the vertex $e$ is exactly $i$. The distance partition of $Q_n$ can be obtained by taking the Hasse diagram of the poset of subsets of an $n$-element set and rotating it clockwise by $90$ degrees;   see Figure~\ref{fig:dist:partition:Q3}.  Each vertex in $X_2(e)$ has a unique set of neighbors in $X_1(e)$.  Hence, since $g$ fixes $X_1(e)$ pointwise and $g$ is an automorphism, $g$ also fixes $X_2(e)$ pointwise.  More generally, for $i \ge 1$, each vertex in $X_{i+1}(e)$ has a unique set of neighbors in $X_i(e)$. It follows by induction on $i$ that $g$ fixes all the vertices of $X$.  Hence, $g$ is the trivial automorphism of $X$.  This proves that $L_e=1$.  It follows that $|G| = |V(X)|~|G_e| = 2^n~ |G_e|$ $\le 2^n~|L_e|~|\Sym(S)| = 2^n (1) (n!)$. 
\qed

A $k$-arc of a graph $X$ is a sequence of $k+1$ vertices $(x_0,x_1,\ldots,x_k)$ such that $x_i$ is adjacent to $x_{i-1}$ for each $i=1,\ldots,k$ and $x_{i+1} \ne x_{i-1}$ for each $i=1,\ldots,k-1$.  A graph $X$ is said to be {\em $k$-arc-transitive} if $\Aut(X)$ acts transitively on the set of $k$-arcs of $X$. A graph $X$ is said to be {\em $k$-arc-regular} if $\Aut(X)$ acts regularly on the set of $k$-arcs of $X$.  A 1-arc-transitive graph is also called {\em arc-transitive}.  Note that an arc-transitive graph is edge-transitive, but the converse doesn't hold (consider $K_{3,4}$).   The hypercube $Q_3$ is 2-arc-transitive but not 3-arc-transitive because it has two kinds of 3-arcs $(x_0,x_1,x_2,x_3)$ - those for which $x_0$ and $x_3$ are adjacent, and those for which $x_0$ and $x_3$ are a distance 3 apart.  Tutte \cite{Tutte:1947} proved the deep result that if $X$ is a cubic, connected, arc-transitive graph, then $X$ is $s$-arc-regular for some $s \le 5$. 

A graph $X=(V,E)$ is said to be {\em distance-transitive} if for any two ordered pairs of vertices $(u,v), (x,y) \in V \times V$,  if the distance between $u$ and $v$ is equal to the distance between $x$ and $y$, then there exists an automorphism of $X$ which maps $(u,v)$ to $(x,y)$.   Let $G:=\Aut(X)$ and fix $v \in V$. Let $X_i(v)$ denote the set of vertices of $X$ whose distance to vertex $v$ is exactly $i$.  It can be shown that a connected graph $X$ of diameter $d$ is distance-transitive if and only if $X$ is vertex-transitive and the stabilizer $G_v$ is transitive on $X_i(v)$ for each $i=1,\ldots,d$.  It can be proved that

\begin{Lemma}
 The hypercube graph $Q_n$ is distance-transitive.
\end{Lemma}

The proof of Theorem~\ref{thm:Qn:autgroup} shows that if $G = \Aut(Q_n) = \Aut(\Cay(\mathbb{Z}_2^n,S))$, then $G_e = \Aut(\mathbb{Z}_2^n,S)$.  Thus, the hypercube $Q_n$ is an example of a Cayley graph $\Cay(\mathbb{Z}_2^n,S)$ that has the smallest possible full automorphism group $R(\mathbb{Z}_2^n) \rtimes \Aut(\mathbb{Z}_2^n,S)$.  In other words, the hypercube is a normal Cayley graph.   

\subsection{Recent results}

In this section, we present some recent results in the literature on the automorphism groups of some families of Cayley graphs.  

The alternating group graph $AG_n$ ($n \ge 3)$ is defined to be the Cayley graph of the alternating group $A_n$ with respect to the generating set $S := \{s_i, s_i^{-1}: i = 3,4,\ldots,n\}$, where $s_i := (1,2,i)$.  Thus, this graph is regular with valency $2(n-2)$.  Zhou \cite{Zhou:2011} obtained the automorphism group of $AG_n$. 

\begin{Theorem}
The automorphism group of the alternating group graph $AG_n=\Cay(A_n,S)$ is equal to $R(A_n) \rtimes \Aut(A_n,S)$, where $\Aut(A_n,S) \cong S_{n-2} \times S_2$. 
\end{Theorem}

Recall that $\mathbb{Z}_2^n$ is the abelian group consisting of vectors of length $n$ over the binary field, with the group operation defined as vector addition (componentwise, mod 2).  

\begin{Definition}
The {\em folded hypercube graph $FQ_n$} is defined to be the Cayley graph of $\mathbb{Z}_2^n$ with respect to the generating set $S = \{e_1,\ldots,e_n,u\}$, where $u= e_1+\ldots+e_n$. 
\end{Definition}

The folded hypercube graph $FQ_n$ can be obtained by starting with the hypercube graph $Q_n$ and adding edges (corresponding to the generator $u$) between diametrically opposite vertices.  Mirafzal  \cite{Mirafzal:2011} obtained the automorphism group of $FQ_n$.  

\begin{Theorem} \label{thm:FQn:autgroup} 
 The automorphism group of the folded hypercube $FQ_n$ ($n \ge 4$) is isomorphic to $\mathbb{Z}_2^n \rtimes S_{n+1}$.  The folded hypercube $FQ_n$ ($n \ge 4$) is a normal Cayley graph and is edge-transitive. 
\end{Theorem}

\begin{Definition}
The {\em augmented cube graph $AQ_n$} ($n \ge 4$) is defined to be the Cayley graph $\Cay(\mathbb{Z}_2^n,S)$, where 
\[S = \{e_1,\ldots,e_n\} \cup \{00 \cdots 00011, 00 \cdots 00111, 00 \cdots 01111, \cdots,  11 \cdots 1111\}. \]
\end{Definition}

Thus, $AQ_n$ is regular with valency $2n-1$.  The automorphism group of $AQ_n$ was investigated in \cite{Choudum:Sunitha:2008} \cite{Ganesan:AutAQn:TechRep:2015}.  It can be shown that $AQ_n$ is a normal Cayley graph and is not edge-transitive.  

We have recalled some symmetry properties of the hypercube $Q_n$, the folded hypercube $FQ_n$ and the augmented cube $AQ_n$, all of which are Cayley graphs of $\mathbb{Z}_2^n$. We pose the following new problems:

\begin{Problem}
Obtain the automorphism group of various other families of Cayley graph of $\mathbb{Z}_2^n$. 
\end{Problem}

\begin{Problem}
Obtain necessary or sufficient conditions on the subset $S \subseteq \mathbb{Z}_2^n$ for the Cayley graph $\Cay(\mathbb{Z}_2^n,S)$ to be edge-transitive.
\end{Problem}

The balanced hypercubes are a particular family of Cayley graphs. It was proved recently that the balanced hypercubes are non-normal Cayley graphs (cf. \cite{Zhou:etal:2015}).  An open problem is to determine the automorphism group of the balanced hypercubes. For recent results on the automorphism groups of Cayley graphs generated by transpositions, see Section~\ref{sec:cayley:transp}. 

\section{Connectivity of transitive graphs}
\label{sec:connec:of:transitive}

In this section, we define the vertex-connectivity and edge-connectivity of graphs, explain their importance in the design of fault-tolerant communication networks, and prove Whitney's theorem that the vertex-connectivity of a graph is at most its edge-connectivity. We use Watkins' theory of atomic parts to prove that an edge-transitive graph has vertex-connectivity equal to its minimum degree.  We recall some further results that were proved using the theory of atomic parts. 

\subsection{Vertex-connectivity and fault-tolerance}

In this section, we define the vertex-connectivity, fault-tolerance, and edge-connectivity of graphs and prove Whitney's theorem that the vertex-connectivity of a graph is at most its edge-connectivity.  We also discuss Menger's theorem and the importance of having parallel paths between any two nodes in a communication network. 

Let $X=(V,E)$ be a simple, undirected graph.  The vertex-connectivity of a graph $X$, denoted by $\kappa(X)$, is the minimal number of vertices of $X$ whose removal leaves a disconnected graph or an isolated vertex. Thus, the vertex-connectivity of the complete graph $K_{n}$ is $n-1$.  The vertex-connectivity of a disconnected graph is defined to be 0. Since removing the set of all neighbors of a vertex isolates the vertex, the vertex-connectivity of a graph is at most the minimum degree $\delta(X)$ of a vertex.  Hence, $\kappa(X) \le \delta(X)$.  A graph is said to be {\em $k$-connected} if its vertex-connectivity is at least $k$. 

In many applications, such as the design of communication networks, it is desired that the network remains connected even if some of the nodes or links in the network fail.  
The fault-tolerance of a graph $X$, denoted by $f(X)$,  is the maximum number of faults (node failures) that can be tolerated without disconnecting the graph (or without leaving a single vertex).    In the definition of this graph invariant, it is assumed that the faulty nodes are chosen by an adversary (this is the worst case scenario).  For example, consider the star graph $K_{1,r}$ ($r \ge 2$).   If a leaf node fails (i.e. is removed from the graph), then the rest of the graph is still connected.  However, if the center node fails, then the network becomes disconnected.  Thus, the star graph cannot tolerate even a single node fault and hence its fault-tolerance is equal to 0.  The fault-tolerance of the cycle graph is equal to 1.  

A set of vertices of a graph $X$ whose removal disconnects $X$ is called a separating set for $X$.  If $W$ is a minimum separating set for $X$, then the graph can tolerate up to $|W|-1$ faults but cannot tolerate $|W|$ faults, and so its fault-tolerance is equal to $|W|-1$.  It is clear that $f(X) = \kappa(X) - 1$, and so the problems of obtaining the vertex-connectivity and fault-tolerance of a graph are equivalent. 

Besides node failures in a communication network, link failures can also occur.  The edge-connectivity of a graph $X$, denoted by $\lambda(X)$, is defined to be the minimal number of edges whose removal disconnects the graph.   A {\em disconnecting set} of $X$ is a set of edges whose removal disconnects the graph.  It is clear that removing the set of all edges incident to a vertex disconnects the vertex from the rest of the graph, whence $\lambda(X) \le \delta(X)$.   Observe that if $F$ is a minimum disconnecting set, then $F$ is of the form $E(S, \overline{S}) := \{ xy: x \in S, y \in \overline{S}, xy \in E(X) \}$ for some subset $S \subseteq V$.  This implies that a minimum disconnecting set for the complete graph $K_n$ has $n-1$ edges. 

Whitney \cite{Whitney:1932} proved that the vertex-connectivity of a graph is at most its edge-connectivity: 

\begin{Theorem} \label{thm:whitney:connectivity}
 Let $X$ be a nontrivial graph. Then
 \[ \kappa(X) \le \lambda(X) \le \delta(X). \]
\end{Theorem}

\begin{figure}
\begin{centering}
\begin{tikzpicture}[scale=1,auto]
\vertex[fill] (x1) at (0,0) [label=left:$x_1$] {};
\vertex[fill] (x2) at (0,-1)  {};
\vertex[fill] (x3) at (0,-2) {};
\vertex[fill] (y1) at (5,0) {};
\vertex[fill] (y2) at (5,-1) {};
\vertex[fill] (y3) at (5,-2) {};

\vertex[fill] (z1) at (5.7,-0.5) {};
\vertex[fill] (z2) at (6.1,-1.5) {};

\path
    (x1) edge (y1) 
    (x1) edge (y2)
    (x2) edge (y2)
    (x3) edge (y3)
;
\draw (0,-1) ellipse (1cm and 2cm);
\draw (5.5,-1) ellipse (1cm and 2cm);
\node at (2.5, 0.5) {$F$};
\node at (0, -3.5) {$S$};
\node at (5.5, -3.5) {$\overline{S}$};
\end{tikzpicture}
\caption{Proof of Whitney's inequality.} 
\label{fig:proof:kappa:le:lambda}
\end{centering}
\end{figure}
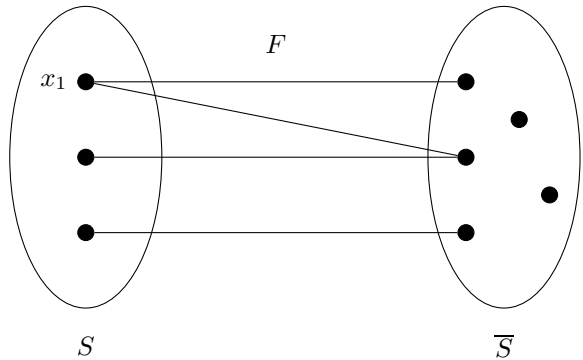

\noindent \emph{Proof:} Let $k = \lambda(X)$ and let $F = \{x_1 y_1,\ldots, x_k y_k\}$ be a set of $k$ edges whose removal disconnects $X$. Here, the $x_i$'s are not necessarily distinct and $F=E(S, \overline{S})$ for some $S$ containing $\{x_1,\ldots,x_k\}$; see  Figure~\ref{fig:proof:kappa:le:lambda}.   If $G - \{x_1,\ldots,x_k\}$ is disconnected, then $\kappa(X) \le k$  and so we are done.  Suppose $G-\{x_1,\ldots,x_k\}$ is connected.  Then, $S$ has no other vertices besides $\{x_1,\ldots,x_k\}$.  In particular, every vertex in $S-\{x_1\}$ is incident to some edge in $F$.  Thus, the degree of $x_1$, which equals the sum of the number of its neighbors in $S$ and the number of its neighbors in $\overline{S}$, is at most $|F|$.  Thus, $d(x_1) \le k$. It follows that $\kappa(X) \le k$. 
\qed

Two $s$-$t$ paths in a graph $X$ are said to be {\em independent} if they have no vertices in common except the end vertices $s$ and $t$.  Let $s$ and $t$ be two nonadjacent vertices of $X$.  If there exists a set of $k$ independent $s$-$t$ paths in $X$, then the number of vertices that need to be removed from $X$ to separate $s$ from $t$ is at least $k$ since at least one vertex needs to be removed from each independent path. Thus, the minimum size of an $s$-$t$ separating set is greater than or equal to the maximum number of independent $s$-$t$ paths.  Menger's theorem \cite{Menger:1927} asserts that this inequality is in fact an equality.  Perhaps the most cited result in graph theory and the result most used in applications of graph theory is Menger's theorem, which is the following.

\begin{Theorem}
 Let $s$ and $t$ be nonadjacent vertices of $X$.  Then, the minimum number of vertices separating $s$ from $t$ is equal to the maximum number of independent $s$-$t$ paths.
\end{Theorem}

Many proofs of Menger's theorem are known (cf. \cite{Bollobas:1998}). The above result is sometimes referred to as the local version of Menger's theorem. The global version of Menger's theorem is the following: the vertex-connectivity of a graph is at least $k$ if and only if every two vertices in the graph are joined by $k$ independent paths. Recall that the vertex-connectivity of a graph was defined in terms of the minimum size of a separating set.  Due to Menger's theorem, one can give a different but equivalent definition of vertex-connectivity in terms of the number of independent paths.  There is also an edge form of Menger's theorem which relates the maximum number of edge-disjoint (rather than node-disjoint) parallel paths to the edge-connectivity of the graph (cf. \cite{Bollobas:1998}). 

If a graph has high fault-tolerance (i.e. $\kappa(X)$ is large), then by Menger's theorem it has a large number of parallel (i.e. independent) paths between any two vertices.  Parallel paths allow for efficient communication since data can be transmitted  simultaneously using parallel paths. Parallel paths also ensure fault-tolerant communication because if the nodes on some of the paths fail, other paths are still available for communication.  

Let $X$ be a graph on vertex set $\{v_1,\ldots,v_n\}$ and with connectivity $\kappa(X)$. Then $|E(X)| = \sum d(v_i) /2 \ge n \delta(X) /2 \ge n \kappa(X) /2$.  Thus, a $k$-connected graph must have at least $\lceil nk/2 \rceil$ edges.  The so-called Harary graphs $H_{k,n}$ attain this bound, i.e. they are examples of $k$-connected graphs that have the minimum possible number of edges $\lceil nk/2 \rceil$ (cf. \cite{Harary:1962}).  If $1 < k < n$ and $k$ is even, the Harary graph $H_{k,n}$ can be constructed as follows: place $n$ vertices on a circle, and join each vertex $x$ to the $k/2$ closest vertices on each side of $x$.  It can be shown that $\kappa(H_{k,n}) = k$.  Observe that the graph constructed is a circulant.  

\subsection{Connectivity of transitive graphs is maximum possible}

Recall that for every nontrivial graph $X$, the vertex-connectivity $\kappa(X)$ and edge-connectivity $\lambda(X)$ satisfy the inequality $\kappa(X) \le \lambda(X) \le \delta(X)$.  A graph for which $\kappa(X) = \delta(X)$ is said to have vertex-connectivity that is maximum possible (or is said to have optimal fault-tolerance). A graph is said to have edge-connectivity that is maximum possible if $\lambda(X) = \delta(X)$.  When we choose a particular graph as the topology of a communication network, it is desirable that the (vertex- or edge-) connectivity of the graph be maximum possible.   Characterizing the graphs for which the vertex-connectivity or edge-connectivity is maximum possible is in general a difficult, interesting and open problem.  

In this section, we show that edge-transitive graphs and vertex-transitive graphs are good topologies to use for communication networks in part because their symmetry properties imply that their connectivity is maximum possible.  We shall use Watkins' theory of atomic parts to prove that the vertex-connectivity of an edge-transitive graph is maximum possible. We state some other known results on optimal fault-tolerance of some families of graphs which were proved using Watkins' theory of atomic parts.

Mader \cite{Mader:1971} proved the following result:

\begin{Theorem}
 If $X$ is a connected vertex-transitive graph, then its edge-connectivity $\lambda(X)$ is equal to its minimum degree $\delta(X)$.
\end{Theorem}

This result settles the question of edge-connectivity for all vertex-transitive graphs and in particular for all Cayley graphs.  Henceforth, we focus our study on the vertex-connectivity of graphs.   

Let $X$ be a simple, undirected graph.  The connectivity of the complete graph is known, and so we assume that the graph $X$ is not the complete graph.  Recall that a separating set for $X$ is a subset of vertices whose removal disconnects the graph.  By definition, a minimum separating set has cardinality $\kappa(X)$.  If $W$ is a minimum separating set for $X$, then the connected components of $X - W$ are called {\em parts} of $X$.  If the subgraph $A$ is a part of $X$, then the set $N(A)$ of vertices in $V(X)-A$ which are adjacent to some vertex in $A$ is a minimum separating set for $X$; we say $N(A)$ is the minimum separating set corresponding to the part $A$. A part is {\em atomic} if it has the smallest possible number of vertices.  Let $p(X)$ denote the number of vertices in an atomic part of $X$.  A graph $X$ is said to be {\em hypo-connected} if $\kappa(X) < \delta(X)$.  An equivalent condition for a graph to be hypo-connected is that each atomic part has at least 2 vertices:

\begin{Lemma} \label{lem:TFAE:atomicparts}
 Let $X$ be a connected graph.  The following are equivalent:
 \\(a) $p(X) \ge 2$
 \\(b) $\kappa(X) < \delta(X)$.
 \\(c) For each vertex $v \in V(X)$, the set of neighbors of the vertex $v$ is not a minimum separating set for $X$. 
\end{Lemma}

\noindent \emph{Proof:} 
\\(a) $\implies$ (b):  Suppose $\kappa(X) = \delta(X)$.  Let $x$ be a vertex having minimal degree $\delta(X)$.  Then, the set of  neighbors of $x$ is a minimum separating set for $X$, whence $\{x\}$ is the vertex set of a part of $X$. Hence $p(X)=1$.
\\(b) $\implies$ (a): Suppose $p(X)=1$ and suppose $\{x\}$ is the vertex set of a part of $X$.  The set of neighbors of $X$ is a minimum separating set for $X$ and has cardinality at least $\delta(X)$.  Thus, $\kappa(X) = \delta(X)$.  
\\The equivalence of (a) and (c) is left to the reader.
\qed

\begin{Theorem} \label{thm:distinct:atomic:parts:disjoint} (Watkins \cite{Watkins:1970}) 
 Let $X$ be a connected graph. Then, the distinct atomic parts of $X$ are vertex-disjoint. 
\end{Theorem}

We recall some definitions from permutation-group theory \cite{Dixon:Mortimer:1996}. Suppose $G$ acts on a set $V$.  A subset $\Delta \subseteq V$ is a {\em block} for $G$ if for each $g \in G$, $\Delta=\Delta^g$ or $\Delta \cap \Delta^g = \phi$.   If $G$ acts transitively on $V$ and $\Delta$ is a block for $G$, then the set $\{\Delta^g: g \in G \}$ of all translates of a block $\Delta$ is a partition of $V$ and is called a {\em complete block system for $G$}.   Theorem~\ref{thm:distinct:atomic:parts:disjoint} implies that:

\begin{Corollary}
 \label{cor:atoms:complete:block:system:AutX}
If $X$ is a connected vertex-transitive graph, then the set of atomic parts of $X$ forms a complete block system for $\Aut(X)$. 
\end{Corollary}

Watkins \cite{Watkins:1970} obtained the following sufficient condition for a graph to have optimal fault-tolerance:  

\begin{Theorem} \label{thm:Watkins:edgetrans:implies:keqd} 
 If $X$ is a connected edge-transitive graph, then its vertex-connectivity $\kappa(X)$ is equal to its minimum degree $\delta(X)$. 
\end{Theorem}

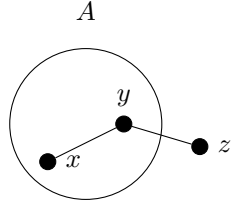
\begin{figure}
\begin{centering}
\begin{tikzpicture}[scale=1,auto]
\vertex[fill] (vx) at (0,0) [label=right:$x$] {};
\vertex[fill] (vy) at (1,0.5) [label=above:$y$] {};
\vertex[fill] (vz) at (2,0.2) [label=right:$z$] {};

\path
    (vx) edge (vy) 
    (vy) edge (vz)
;
\draw (0.5,0.5) ellipse (1cm and 1cm);
\node at (0.5, 2) {$A$};
\end{tikzpicture}
\caption{Proof of Theorem~\ref{thm:Watkins:edgetrans:implies:keqd}.} 
\label{fig:proof:edgetrans:implies:keqd}
\end{centering}
\end{figure}

\noindent \emph{Proof:} Let $X$ be a connected edge-transitive graph.  By way of contradiction, suppose $\kappa(X) < \delta(X)$. By Lemma~\ref{lem:TFAE:atomicparts}, $p(X) \ge 2$.  Let $A$ be an atomic part of $X$.  Then $A$ is an induced subgraph of $X$ on 2 or more vertices and is connected.  Also, there exists a minimum separating set $W$ for $X$ such that $A$ is a connected component of $X-W$.  Since $X$ and $A$ are each connected, there exist vertices $x, y, z$ such that $x, y \in V(A)$, $z \in V(X)-V(A)$, and $xy, yz \in E(X)$ (see Figure~\ref{fig:proof:edgetrans:implies:keqd}).  By the edge-transitivity of $X$, there is an automorphism $\phi$ of $X$ that takes edge $xy$ to edge $yz$.  Since $\phi$ is an automorphism of $X$, $A^\phi$ is also an atomic part of $X$.  But $A^\phi$ has partial overlap with $A$, contradicting Theorem~\ref{thm:distinct:atomic:parts:disjoint}.  Thus, $p(X)=1$, which implies $\kappa(X) = \delta(X)$. 
\qed

Watkins \cite{Watkins:1970} and Mader \cite{Mader:1970} independently obtained a lower bound on the vertex-connectivity of a vertex-transitive graph in terms of the valency  $\delta(X)$:

\begin{Theorem}
If $X$ is a connected vertex-transitive graph, then its vertex-connectivity $\kappa(X)$ is at least $\frac{2}{3} (\delta(X)+1)$. 
\end{Theorem}

Thus, the vertex-connectivity of every connected vertex-transitive graph is bounded as
\[
 \frac{2}{3} (\delta(X)+1) \le \kappa(X) \le \lambda(X) \le \delta(X).
\]

Another sufficient condition for optimal vertex-connectivity was obtained by Mader \cite{Mader:1970}:

\begin{Theorem} \label{thm:Mader:K4:free}
 If $X$ is a connected vertex-transitive graph which does not contain a $K_4$, then its vertex-connectivity $\kappa(X)$ is equal to its minimum degree $\delta(X)$. 
\end{Theorem}

Sufficient conditions for a graph to have vertex-connectivity equal to the minimum degree are given in Theorem~\ref{thm:Watkins:edgetrans:implies:keqd} and Theorem~\ref{thm:Mader:K4:free}.  However, there exist graphs which do not satisfy the hypotheses of these assertions and which still have vertex-connectivity that is maximum possible; for example, some families of circulants and the family of augmented cubes are neither edge-transitive nor $K_4$-free but have vertex-connectivity that is maximum possible (see Section~\ref{sec:kappaeqdelta}).   Obtaining necessary and sufficient conditions for a particular family of transitive graphs to have vertex-connectivity that is maximum possible is in general an interesting and open problem.   Boesch and Tindell \cite{Boesch:Tindell:1984} characterized the circulants which have vertex-connectivity that is maximum possible. In the next section, we mention other results and some new open questions.  

Imrich \cite{Imrich:1979} weakened the condition of edge-transitivity in Theorem~\ref{thm:Watkins:edgetrans:implies:keqd}: 

\begin{Theorem} \label{thm:Imrich:connec:1979}
 Let $X$ be a connected graph. For $f \in E(X)$, let $X_f$ be the component of the edge orbit $f^{\Aut(X)}$ which contains the edge $f$.  Suppose that for every edge $f \in E(X)$,
 \\(i) The edge orbit $f^{\Aut(X)}$ is a spanning subgraph of $X$, and 
 \\(ii) $\delta(X) \le 2 |V(X_f)|$ or $|V(X)| \le 3 |V(X_f)|$.
 \\Then, $\kappa(X) = \delta(X)$. 
\end{Theorem}

Imrich \cite{Imrich:1979} used this result to prove the optimal fault-tolerance of some Cayley graphs of the symmetric group:

\begin{Theorem} \label{thm:Imrich:cay:conjugacy:kappa:eq:delta}
 Let $n \ge 5$.  Let $S$ be the union of some conjugacy classes $C_1, \ldots, C_r$ of $S_n$ such that $S$ contains an odd permutation.  Then, the vertex-connectivity of the Cayley graph $\Cay(S_n,S)$ is maximum possible. 
\end{Theorem}

A special case of Theorem~\ref{thm:Imrich:cay:conjugacy:kappa:eq:delta} occurs when $r=1$ and $C_1$ is the set of all ${n \choose 2}$ transpositions in $S_n$.  In this case, the Cayley graph $\Cay(S_n,S)$ is called the complete transposition graph, and by Theorem~\ref{thm:Imrich:cay:conjugacy:kappa:eq:delta} its vertex-connectivity is equal to ${n \choose 2}$.  But in this case, the Cayley graph $\Cay(S_n,S)$ is also edge-transitive. An open question is to characterize the families of sets $S$ for which $\Cay(S_n,S)$ is edge-transitive:

\begin{Problem}
 Let $S$ be the union of conjugacy classes $C_1,\ldots,C_r$ of the symmetric group $S_n$ such that $S$ contains an odd permutation.  
 \\(a) Characterize the families of sets $S$ for which $\Cay(S_n,S)$ is edge-transitive.
 \\(b) For various families $S$, determine the automorphism group of the Cayley graph $\Cay(S_n,S)$. 
\end{Problem}

In the special case when $S$ is the set of all cyclic permutations in $S_n$, $S$ is the union of exactly $n-1$ conjugacy classes of $S_n$.  In this case, $\Cay(S_n,S)$ is called the assignment polytope, and its vertex-connectivity is maximum possible by Theorem~\ref{thm:Imrich:cay:conjugacy:kappa:eq:delta}. 

Let $H$ be a group and let $S \subseteq H$.  Suppose $1 \notin S = S^{-1}$ and $S$ generates $H$.  We say $S$ is a {\em minimal generating set for $H$} if $\langle S - \{h,h^{-1} \} \rangle$ is a proper subgroup of $H$ for each $h \in S$. Godsil \cite{Godsil:connec:minimalCayley:1981} proved the following:

\begin{Theorem} \label{thm:Godsil:minimalgenset:connec}
 Let $H$ be a group.  If $S$ is a minimal generating set for $H$, then the vertex-connectivity of the Cayley graph $\Cay(H,S)$ is maximum possible.
\end{Theorem}

It follows that if $S$ is a set of transpositions generating $S_n$ such that the transposition graph $T(S)$ is a tree, then the vertex-connectivity of $\Cay(S_n,S)$ is equal to $n-1$. The optimal fault-tolerance of Cayley graphs whose generator sets satisfy other properties (such as hierarhical or quasi-minimal) is proved in \cite{Akers:Krishnamurthy:1987} \cite{Alspach:1992}.  

The results mentioned so far on the connectivity of graphs  pertain to either the vertex-connectivity or the edge-connectivity of graphs.   One can define various other notions of connectivity or fault-tolerance, and some of these definitions are motivated by applications.  If $P$ is any graph-theoretic property, the {\em conditional connectivity of $X$ with respect to $P$}, denoted by $\kappa(X; P)$ is the smallest size of a subset $W \subseteq V(X)$ such that every component of the disconnected graph $X-W$ has property $P$.  For example, $P$ can be the property that the graph has a cycle, in which case $\kappa(X;P)$ is called the {\em cyclic connectivity of $X$} and is equal to the smallest number of vertices whose removal disconnects $X$ and such that every component of the disconnected graph contains a cycle.  For an overview of this topic of conditional connectivity, see Harary \cite{Harary:1983}. In a recent paper \cite{Zhou:etal:2015}, various reliability or conditional connectivity measures of the balanced hypercube are obtained.
\section{Automorphisms and optimal fault-tolerance of some families of graphs}
\label{sec:kappaeqdelta}

In this section, we prove that the vertex-connectivity of many families of Cayley graphs is maximum possible.  We consider the family of hypercube graphs, the family of folded hypercubes, Cayley graphs generated by transpositions and Cayley graphs from linear codes.  

\subsection{Connectivity of the hypercube}

In this section, we give three proofs that the vertex-connectivity of the hypercube is maximum possible. The first proof uses information about the automorphism group of the graph and is non-constructive;  the second proof is constructive and the third proof is by induction.

The hypercube $Q_n$ is defined to be the graph on vertex set $\{0,1\}^n$, and two binary strings $x=x_1 \cdots x_n$ and $y=y_1 \cdots y_n$ are adjacent vertices in $Q_n$ if and only if they differ in exactly one coordinate.  There are other equivalent definitions of the hypercube.  Recall that the hypercube is isomorphic to the Cayley graph of the permutation group generated by $n$ disjoint transpositions (cf. Example~\ref{eg:hypercube:transpositions}), and so the hypercube graph could have also been defined as a particular kind of Cayley graph.  

Let $\mathbb{F}_2^n$ be the $n$-dimensional vector space over the binary field $\mathbb{F}_2$.   The set of unit vectors $e_i$ ($i=1,\ldots,n)$ is a basis for the vector space $\mathbb{F}_2^n$.    Recall that $\mathbb{F}_2^n$ is an abelian group $\mathbb{Z}_2^n$ under the operation of vector addition, and the subgroups of $\mathbb{Z}_2^n$ correspond to the subspaces of the vector space.  Note that the Cayley graph of the abelian group $\mathbb{Z}_2^n$ with respect to the set of $n$ unit vectors $e_1,\ldots,e_n$ is isomorphic to the hypercube graph $Q_n$. In the sequel, we view $Q_n$ as the Cayley graph $\Cay(\mathbb{Z}_2^n,S)$, where $S = \{e_1,\ldots,e_n\}$. 

The hypercube $Q_n$ is an $n$-regular, vertex-transitive graph on $2^n$ vertices.  For $x, y \in V(Q_n) = \mathbb{Z}_2^n$, $xy$ is an edge of $Q_n$ iff $x+y = e_i$ for some $i$; this edge is said to have edge label (or color) $e_i$ or to be of dimension $i$.  If $y = 1 \cdots 1 0 \cdots 0$ is a vertex consisting of $k$ 1's and $n-k$ 0's, then the distance from $y$ to the identity vertex $e = 0 \cdots 0$ is exactly $k$.  The path from $e$ to $y$ can be described by a sequence $(e_1, e_2, \ldots, e_k)$ of labels of the edges on the path.  

Recall that for any graph $X$, the vertex-connectivity $\kappa(X)$ is at most the minimal degree $\delta(X)$, and graphs for which $\kappa(X)$ is equal to $\delta(X)$ are said to have optimal fault-tolerance (or to have vertex-connectivity that is maximum possible).  We now establish that the hypercube graphs have optimal fault-tolerance.  We shall give three proofs of this result.   The first proof establishes that $Q_n$ is edge-transitive and then appeals to Watkins' theorem that the vertex-connectivity of an edge-transitive graph is maximum possible.  This proof is non-constructive in the sense that it establishes only the existence of $n$ independent paths between any two vertices.  In the second proof, we give a construction for a set of $n$ independent paths between any two vertices in $Q_n$, thereby establishing that $\kappa(Q_n) \ge n$.  The third proof is by induction. 

\begin{Theorem} \label{thm:Qn:kappa:eq:delta}
 The vertex-connectivity of the hypercube graph $Q_n$ is equal to $n$.
\end{Theorem}

\noindent \emph{First proof:} 
Recall Theorem~\ref{thm:Watkins:edgetrans:implies:keqd} that a sufficient condition for the vertex-connectivity of a graph to equal its minimum degree is that the graph be edge-transitive. 
Hence, it suffices to prove that the hypercube $Q_n$ is edge-transitive.   

Let $G$ be the automorphism group of the hypercube $Q_n$.  Let $\{x,x+e_j\}$ be an edge of $Q_n$.  It suffices to show that there is an element in $G$ that maps the edge $\{e,e_1\}$ to the edge $\{x, x+e_j\}$.  There exists an element $g \in G$ which permutes the unit (basis) vectors and which takes the edge $\{e, e_1\}$ to the edge $\{e,e_j\}$.   The translation map $r_x$ which takes each vertex $u$ to $u+x$ is also an automorphism of the hypercube graph $Q_n$.  The composition $g r_x$ takes the edge $\{e,e_1\}$ to the edge $\{x,x+e_j\}$.  Hence, the action of $G$ on $E(Q_n)$ has a single orbit. 
\qed

\bigskip
\noindent \emph{Second proof of Theorem~\ref{thm:Qn:kappa:eq:delta}:}
It suffices to show that $\kappa(Q_n) \ge n$, for the opposite inequality holds trivially.  Let $x$ and $y$ be any two distinct vertices of $Q_n$.  It suffices to show that there exist $n$ independent $x$-$y$ paths in the graph $Q_n$.  We can assume without loss of generality that $x$ is the identity vertex $e$ because the translations $r_z$ that take each vertex $u$ to $u+z$ are automorphisms of $X$.  We can also assume that $y$ is a vector of the form $1 \cdots 1 0 \cdots 0$ where all the 1's precede all the 0's because permutations of the set of unit vectors $e_i$ ($i=1,\ldots,n$) induce automorphisms of the graph $Q_n$ that fix the identity vertex $e$.  Suppose the number of 1's in $y$ is equal to $k$.  Then, a shortest path from $e$ to $y$ is the path (of length $k$) corresponding to the sequence of edge labels $(e_1,e_2,\ldots,e_k)$.   A left cyclic shift of this sequence is the sequence $(e_2,e_3,\ldots,e_k,e_1)$, which corresponds to another path from $e$ to $y$. The set of all cyclic shifts of the sequence $(e_1,\ldots,e_k)$ gives a set $A$ of $k$ independent paths from $e$ to $y$. The set of sequences $(e_i, e_1, e_2, \ldots, e_k, e_i)$ $(i=k+1,\ldots,n)$ corresponds to a set $B$ of $n-k$ independent paths from $e$ to $y$.  The union $A \cup B$ is a set of $n$ independent paths from $e$ to $y$.
\qed

\bigskip
\noindent \emph{Third proof of Theorem~\ref{thm:Qn:kappa:eq:delta}:}
We prove that $\kappa(Q_n)=n$ by induction on $n$.  If $n=1$, then $Q_1$ is the graph $K_2$, which clearly has vertex-connectivity equal to 1.  Let $n \ge 2$ and assume that the assertion holds for smaller values of $n$.  The graph $Q_n$ consists of two disjoint copies $Y, Y'$ of $Q_{n-1}$ which are joined to each other by a perfect matching consisting of $2^{n-1}$ edges.  Let $W$ be a separating set for $Q_n$.  We show that $|W| \ge n$.   We consider two cases.  

First, suppose $Y-W$ and $Y'-W$ are both connected.  Then, the only way $Q_n-W$ can be disconnected is if $W$ contains at least one vertex incident to each edge of the perfect matching. Hence, $|W| \ge 2^{n-1} \ge n$, as was to be shown.

Now suppose at least one of $Y-W$ or $Y'-W$ is disconnected, say $Y-W$ is disconnected. By the inductive hypothesis, $|W \cap V(Y)| \ge n-1$.  If $W$ does not contain any vertex of $Y'$, then $Q_n-W$ will be connected due to the perfect matching, a contradiction.  Hence, $W$ contains at least one vertex of $Y'$ and hence contains at least $n$ vertices in total.
\qed

\subsection{Containers and parallel paths in the folded hypercube}

El-Amawy and Latifi \cite{ElAmawy:Latifi:1991} proposed the folded hypercube graph as a topology for interconnection networks. The folded hypercube graph $FQ_n$ ($n \ge 2$) is defined to be Cayley graph $\Cay(\mathbb{Z}_2^n,S)$, where $\mathbb{Z}_2^n$ is the abelian group consisting of all $0$-$1$ vectors of length $n$ (with mod 2, componentwise addition) and the generating set $S = \{e_1,\ldots,e_n,u\}$, with $u=e_1+\ldots,e_n$.  In other words, the folded hypercube $FQ_n$ is obtained by taking the hypercube $Q_n$ and adding edges (corresponding to the generator $u$) which join each vertex to its diametrically opposite vertex.  The motivation for adding these complementary edges to the hypercube is that they reduce the diameter of the graph from $n$ to about $n/2$.  For if two vertices in $FQ_n$ differ in more than half of the coordinates, a shorter path between these two vertices is obtained by using the complementary edge.  For example, the length of a shortest path in $FQ_6$ from vertex $e=000000$ to vertex $011111$ is 2; one such shortest path is the path corresponding to the sequence of edge labels (or generators) $(u,e_1)$.  The folded hypercube $FQ_n$ is a regular graph of valency $n+1$.  Thus, its vertex-connectivity satisfies $\kappa(FQ_n) \le n+1$.  

\begin{Theorem} \label{thm:FQn:kappa:eq:delta}
 The vertex-connectivity of the folded hypercube $FQ_n$ ($n \ge 4$) is equal to $n+1$. 
\end{Theorem}

\noindent \emph{Proof:}
We construct a set of $n+1$ independent paths between any two vertices $x, y \in V(FQ_n)$.  By vertex-transitivity of $FQ_n$, we can assume that $x$ is the identity vertex $e$.   Suppose the number of 1's in $y$ is $r$.  We can assume without loss of generality that the $r$ 1's precede the $n-r$ 0's because each permutation $g \in S_n$ of the coordinates of the vectors acts on $V(FQ_n)$ as a group of automorphisms of $FQ_n$ (cf. Theorem~\ref{thm:FQn:autgroup}).   We show that there exists a set of $n+1$ independent $e$-$y$ paths in $FQ_n$, where $y$ is the vector $1 \cdots 1 0 \cdots 0$ consisting of $r$ 1's and $n-r$ 0's.   We consider two cases, depending on the value of $r$.

Suppose $r \le \lceil n/2 \rceil$. The path from $e$ to $y$ corresponding to the sequence of edge labels $(e_1,e_2,\ldots,e_r)$ has length $r$.  The set of all cyclic shifts of this sequence gives a set of $r$ independent paths from $e$ to $y$, each of length $r$.  The set of sequences of edge labels $(e_i,e_1,e_2,\ldots,e_r,e_i)$ (for $i=r+1, \ldots, n$) gives a set of $n-r$ independent paths from $e$ to $y$, each of length $r+2$.  Finally, the path corresponding to the sequence $(u, e_1, e_2, \ldots, e_r, u)$ is another path from $e$ to $y$ and has length $r+2$.  It can be verified that these $n+1$ paths are pairwise independent.  

Suppose $r > \lceil n/2 \rceil$.   We can assume without loss of generality that $y$ is a vertex of the form $0 \cdots 0 1 \cdots 1$ consisting of $n-r$ 0's followed by $r$ 1's.  The set of paths corresponding to the sequence of edge labels $(u, e_1, e_2, \ldots, e_{n-r})$ and its cyclic shifts gives a set of $n-r+1$ independent paths from $e$ to $y$, each of length $n-r+1$.  The  set of sequences $(e_i, u, e_1, e_2, \ldots, e_{n-r}, e_i)$ (for $i=n-r+1, \ldots, n$) gives a set of $r$ paths from $e$ to $y$, each of length $n-r+3$.  Altogether, we have a set of $n+1$ independent paths from $e$ to $y$.
\qed

A {\em container} is a set of independent paths in a graph between two given vertices.  The width of the container is the number of independent paths in the container.  The length of the container is the longest length of a path in the container.  Recall from the second proof of Theorem~\ref{thm:Qn:kappa:eq:delta} that we constructed a set of $n$ independent paths in $Q_n$ between any two vertices.  The proof showed that if $x, y \in V(Q_n)$ and $r$ is the distance in $Q_n$ between $x$ and $y$, then there exists a set of $n$ independent paths between $x$ and $y$ such that $r$ of the paths are of length $r$ each and $n-r$ of the paths are of length $r+2$ each.  Thus, we have the following:

\begin{Corollary} 
Given any two vertices in the hypercube graph $Q_n$, we can establish a container of width $n$ and length $r+2$ between these two vertices, where $r$ denotes the number of coordinates where the two vertices differ. 
\end{Corollary}

From the proof of Theorem~\ref{thm:FQn:kappa:eq:delta}, we can deduce the following result.

\begin{Corollary}
 Let $x, y$ be two distinct vertices of the folded hypercube graph $FQ_n$ ($n \ge 4$).  Let $r$ be the number of coordinates where $x$ and $y$ differ.  If $r \le \lceil n/2 \rceil$, then we can establish a container of width $n+1$ and length $r+2$ between vertices $x$ and $y$.  If $r > \lceil n/2 \rceil$, then we can establish a container of width $n+1$ and length at most $\lfloor n/2 \rfloor + 2$ between vertices $x$ and $y$.  
\end{Corollary}

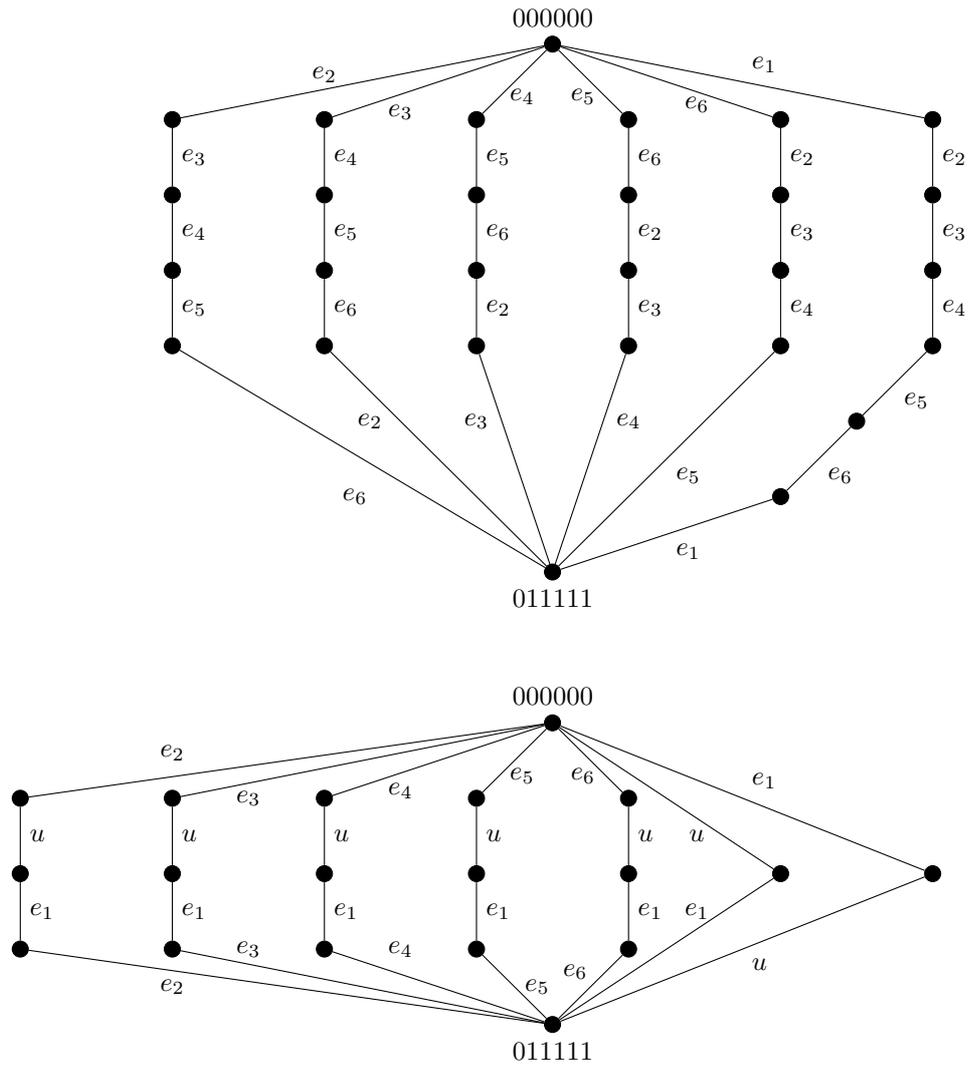
\begin{figure}
\begin{centering}
\begin{tikzpicture}[scale=1,auto]
\vertex[fill] (ve) at (0,0) [label=above:$000000$] {};
\vertex[fill] (vy) at (0,-7) [label=below:$011111$] {};

\vertex[fill] (c1) at (-5,-1) {};
\vertex[fill] (c2) at (-5,-2) {};
\vertex[fill] (c3) at (-5,-3) {};
\vertex[fill] (c4) at (-5,-4) {};
\draw [-] (ve) edge node {} (c1);  \node at (-3,-0.4) {$e_2$};
\draw [-] (c1) edge node {$e_3$} (c2);
\draw [-] (c2) edge node {$e_4$} (c3);
\draw [-] (c3) edge node {$e_5$} (c4);
\draw [-] (c4) edge node {} (vy);  \node at (-2.6,-6) {$e_6$};

\vertex[fill] (a1) at (-3,-1) {};
\vertex[fill] (a2) at (-3,-2) {};
\vertex[fill] (a3) at (-3,-3) {};
\vertex[fill] (a4) at (-3,-4) {};
\draw [-] (ve) edge node {} (a1); \node at (-2,-0.9) {$e_3$};
\draw [-] (a1) edge node {$e_4$} (a2);
\draw [-] (a2) edge node {$e_5$} (a3);
\draw [-] (a3) edge node {$e_6$} (a4);
\draw [-] (a4) edge node {} (vy); \node at (-2.4,-5) {$e_2$};

\vertex[fill] (b1) at (-1,-1) {};
\vertex[fill] (b2) at (-1,-2) {};
\vertex[fill] (b3) at (-1,-3) {};
\vertex[fill] (b4) at (-1,-4) {};
\draw [-] (ve) edge node {} (b1); \node at (-0.4,-0.7) {$e_4$};
\draw [-] (b1) edge node {$e_5$} (b2); 
\draw [-] (b2) edge node {$e_6$} (b3);
\draw [-] (b3) edge node {$e_2$} (b4);
\draw [-] (b4) edge node {} (vy); \node at (-1,-5) {$e_3$};

\vertex[fill] (d1) at (1,-1) {};
\vertex[fill] (d2) at (1,-2) {};
\vertex[fill] (d3) at (1,-3) {};
\vertex[fill] (d4) at (1,-4) {};
\draw [-] (ve) edge node {} (d1); \node at (0.4,-0.7) {$e_5$};
\draw [-] (d1) edge node {$e_6$} (d2);
\draw [-] (d2) edge node {$e_2$} (d3);
\draw [-] (d3) edge node {$e_3$} (d4);
\draw [-] (d4) edge node {} (vy); \node at (1,-5) {$e_4$};

\vertex[fill] (e1) at (3,-1) {};
\vertex[fill] (e2) at (3,-2) {};
\vertex[fill] (e3) at (3,-3) {};
\vertex[fill] (e4) at (3,-4) {};
\draw [-] (ve) edge node {} (e1); \node at (1.9,-0.8) {$e_6$};
\draw [-] (e1) edge node {$e_2$} (e2);
\draw [-] (e2) edge node {$e_3$} (e3);
\draw [-] (e3) edge node {$e_4$} (e4);
\draw [-] (e4) edge node {$e_5$} (vy);

\vertex[fill] (f1) at (5,-1) {};
\vertex[fill] (f2) at (5,-2) {};
\vertex[fill] (f3) at (5,-3) {};
\vertex[fill] (f4) at (5,-4) {};
\vertex[fill] (f5) at (4,-5) {};
\vertex[fill] (f6) at (3,-6) {};
\draw [-] (ve) edge node {$e_1$} (f1);
\draw [-] (f1) edge node {$e_2$} (f2);
\draw [-] (f2) edge node {$e_3$} (f3);
\draw [-] (f3) edge node {$e_4$} (f4);
\draw [-] (f4) edge node {$e_5$} (f5);
\draw [-] (f5) edge node {$e_6$} (f6);
\draw [-] (f6) edge node {$e_1$} (vy);

\pgfmathsetmacro{\phi}{9}
\vertex[fill] (xve) at (0,0-\phi) [label=above:$000000$] {};
\vertex[fill] (xvy) at (0,-4-\phi) [label=below:$011111$] {};

\vertex[fill] (xg1) at (-7,-1-\phi) {};
\vertex[fill] (xg2) at (-7,-2-\phi) {};
\vertex[fill] (xg3) at (-7,-3-\phi) {};
\draw [-] (xve) edge node {} (xg1);  \node at (-5,-0.4-\phi) {$e_2$};
\draw [-] (xg1) edge node {$u$} (xg2);
\draw [-] (xg2) edge node {$e_1$} (xg3);
\draw [-] (xg3) edge node {} (xvy); \node at (-5,-3.5-\phi) {$e_2$};

\vertex[fill] (xc1) at (-5,-1-\phi) {};
\vertex[fill] (xc2) at (-5,-2-\phi) {};
\vertex[fill] (xc3) at (-5,-3-\phi) {};
\draw [-] (xve) edge node {} (xc1);  \node at (-4,-1-\phi) {$e_3$};
\draw [-] (xc1) edge node {$u$} (xc2);
\draw [-] (xc2) edge node {$e_1$} (xc3);
\draw [-] (xc3) edge node {} (xvy);  \node at (-4,-3-\phi) {$e_3$};

\vertex[fill] (xa1) at (-3,-1-\phi) {};
\vertex[fill] (xa2) at (-3,-2-\phi) {};
\vertex[fill] (xa3) at (-3,-3-\phi) {};
\draw [-] (xve) edge node {} (xa1); \node at (-2,-0.9-\phi) {$e_4$};
\draw [-] (xa1) edge node {$u$} (xa2);
\draw [-] (xa2) edge node {$e_1$} (xa3);
\draw [-] (xa3) edge node {} (xvy); \node at (-2,-3-\phi) {$e_4$};

\vertex[fill] (xb1) at (-1,-1-\phi) {};
\vertex[fill] (xb2) at (-1,-2-\phi) {};
\vertex[fill] (xb3) at (-1,-3-\phi) {};
\draw [-] (xve) edge node {} (xb1); \node at (-0.4,-0.7-\phi) {$e_5$};
\draw [-] (xb1) edge node {$u$} (xb2); 
\draw [-] (xb2) edge node {$e_1$} (xb3);
\draw [-] (xb3) edge node {} (xvy); \node at (-0.2,-3.5-\phi) {$e_5$};

\vertex[fill] (xd1) at (1,-1-\phi) {};
\vertex[fill] (xd2) at (1,-2-\phi) {};
\vertex[fill] (xd3) at (1,-3-\phi) {};
\draw [-] (xve) edge node {} (xd1); \node at (0.4,-0.7-\phi) {$e_6$};
\draw [-] (xd1) edge node {$u$} (xd2);
\draw [-] (xd2) edge node {$e_1$} (xd3);
\draw [-] (xd3) edge node {} (xvy); \node at (0.3,-3.3-\phi) {$e_6$};

\vertex[fill] (xe1) at (3,-2-\phi) {};
\draw [-] (xve) edge node {} (xe1); \node at (1.9,-1.5-\phi) {$u$};
\draw [-] (xe1) edge node {} (xvy); \node at (1.9,-2.5-\phi) {$e_1$};

\vertex[fill] (xf1) at (5,-2-\phi) {};
\draw [-] (xve) edge node {$e_1$} (xf1);
\draw [-] (xf1) edge node {$u$} (xvy);

\end{tikzpicture}
\caption{A container between vertices $000000$ and $011111$ in (a) the hypercube $Q_6$, (b) the folded hypercube $FQ_6$.}
\label{fig:container:Q6:FQ6}
\end{centering} 
\end{figure}

Thus, the folded hypercube has wider and shorter containers than the hypercube $Q_n$.  This is illustrated in Figure~\ref{fig:container:Q6:FQ6}(a) and Figure~\ref{fig:container:Q6:FQ6}(b).  Figure~\ref{fig:container:Q6:FQ6}(a) shows a container in the hypercube $Q_6$ between vertices $e=000000$ and $y=0111111$, and Figure~\ref{fig:container:Q6:FQ6}(b) shows a container in the folded hypecube $FQ_6$ between the same two vertices.  In these figures, the edges are labeled by the generators of the corresponding Cayley graph. 

The {\em quality} of a container is defined to be the ratio of the width of the container to the average length of the paths in the container.  For example, if $x, y \in V(Q_n)$ differ in exactly $r$ coordinates, we showed there exists a container containing $r$ paths of length $r$ each and $n-r$ paths of length $r+2$ each.  Thus, there exists a container between $x$ and $y$ of quality $\frac{n}{\{ r \cdot r + (n-r)(r+2)\}/n} = \frac{n^2}{n(r+2)-2r}$.   

The proof of Theorem~\ref{thm:FQn:kappa:eq:delta} on the optimal vertex-connectivity of the folded hypercube $FQ_n$ gave a set of $n+1$ independent paths between any two vertices.  Since an edge-transitive graph has optimal vertex-connectivity, the optimal vertex-connectivity of $FQ_n$ can also be deduced from Theorem~\ref{thm:FQn:autgroup}.  

We recall again the definitions of three families of Cayley graphs of $\mathbb{Z}_2^n$ mentioned above.  The hypercube $Q_n$ is defined to be the Cayley graph $\Cay(\mathbb{Z}_2^n,S)$, where $S = \{e_1,\ldots,e_n\}$.  The folded hypercube $FQ_n$ is the Cayley graph $\Cay(\mathbb{Z}_2^n,S)$, where $S = \{e_1,\ldots,e_n,u\}$ and $u=e_1+\cdots+e_n$.  The augmented cube $AQ_n$ is the Cayley graph $\Cay(\mathbb{Z}_2^n,S)$, where 
\[S = \{e_1,\ldots,e_n\} \cup \{00 \cdots 00011, 00 \cdots 00111, 00 \cdots 01111, \cdots,  11 \cdots 1111\}.\]  

The hypercubes and folded hypercubes are edge-transitive, and so by Theorem~\ref{thm:Watkins:edgetrans:implies:keqd} they have optimal fault-tolerance.  The hypercubes and folded hypercubes are also $K_4$ free, and so by Theorem~\ref{thm:Mader:K4:free} we get another proof of their optimal fault-tolerance.  The augmented cubes have optimal fault-tolerance (cf. \cite[Proposition 4.1]{Choudum:Sunitha:2002}).  However, the augmented cube graphs are neither edge-transitive nor $K_4$-free.  In view of these observations (which are summarized in Table~\ref{table:Qn:FQn:AQn:kappa:eq:delta}), we ask the question of whether there are other sufficient conditions on $S$ (which the augmented cubes would satisfy) for a Cayley graph $\Cay(\mathbb{Z}_2^n,S)$ to have optimal fault-tolerance. 

\begin{table}
\begin{center}
\begin{tabular}{|c|c|c|c|}
\hline
 graph & edge-transitive & $K_4$-free & optimal fault-tolerance  \\
 \hline
 $Q_n$ & yes & yes & yes \\ \hline 
 $FQ_n$ & yes & yes & yes \\ \hline 
 $AQ_n$ & no & no & yes  \\ \hline
  \end{tabular}
\caption{On obtaining new sufficient conditions for optimal fault-tolerance.}
\end{center}
 \end{table}
\label{table:Qn:FQn:AQn:kappa:eq:delta}

\begin{Problem}
 Obtain necessary or sufficient conditions on $S \subseteq \mathbb{Z}_2^n$ for $\Cay(\mathbb{Z}_2^n,S)$ to have optimal fault-tolerance. 
\end{Problem}

\subsection{Automorphism group of Cayley graphs generated by transpositions}
\label{sec:cayley:transp}

In this section, we define the family of Cayley graphs generated by transpositions, mention some recent results on automorphism groups of Cayley graphs generated by transposition, and state some new problems and new open conjectures.

\begin{Definition}
Let $S$ be a set of transpositions in the symmetric group $S_n$.   The {\em transposition graph of $S$}, denoted by $T(S)$, is defined to be the graph with vertex set $\{1,\ldots,n\}$, and two vertices $i$ and $j$ are adjacent in $T(S)$ whenever $(i,j) \in S$.   
\end{Definition}

Thus, the set $S$ of transpositions in $S_n$ can be represented by the (edge set of the) graph $T(S)$ on $n$ vertices.  

\begin{Lemma} \label{lem:GR:transpositiongraphs}
 Let $S$ be a set of transpositions in $S_n$.  Then, 
 \\(a) $S$ generates $S_n$ if and only if the transposition graph $T(S)$ is connected.  
 \\(b) $S$ is a minimal generating set for $S_n$ if and only if the transposition graph $T(S)$ is a tree.
\end{Lemma}

\noindent \emph{Proof:}  Let $S$ be a set of transpositions in $S_n$. We first prove that if the transposition graph $T(S)$ is connected, then $S$ generates $S_n$. Since every permutation in $S_n$ is a product of transpositions, it suffices to show that each transposition in $S_n$ is a product of elements of $S$. Let $(i,j)$ be a transposition in $S_n$.  Since $T(S)$ is connected, there exists a path in $T(S)$ from vertex $i$ to vertex $j$, say $i=x_0, x_1, \ldots, x_r = j$.  The transposition $(i,j)$ is the composition of the following transpositions in $S$: $(x_0,x_1), (x_1,x_2), \ldots, (x_{r-1}, x_r), (x_{r-2}, x_{r-1}), \ldots, (x_0, x_1)$.  To prove the converse, suppose the transposition graph $T(S)$ is not connected.  Suppose $i$ and $j$ are in different connected components of $T(S)$.  Then, the transposition $(i,j) \in S_n$ is not a product of elements of $S$. This proves (a).  A tree is precisely a minimal connected graph, and so (b) follows. 
\qed

Let $S$ be a set of transpositions in $S_n$.  The graph $\Cay(S_n,S)$ is called a {\em Cayley graph generated by transpositions}. If $n$ is even, say $n=2k$, and $T(S)$ is the graph $k K_2$ consisting of $k$ independent edges, then the Cayley graph $\Cay(\langle S \rangle, S)$ is isomorphic to the hypercube graph $Q_n$.  Various families of Cayley graphs generated by transpositions have been well-studied and they have specific names \cite{Lakshmivarahan:etal:1993} \cite{Heydemann:1997}:

\begin{Definition}
 Let $S$ be a set of transpositions in $S_n$. Let $T(S)$ denote the transposition graph of $S$.  
 \\(a) If $T(S)$ is the star $K_{1,n-1}$, then $\Cay(S_n,S)$ is the {\em star graph}.
 \\(b) If $T(S)$ is the path graph $P_n$ on $n$ vertices, then $\Cay(S_n,S)$ is called the {\em bubble-sort graph}.
 \\(c) If $T(S)$ is the cycle graph $C_n$, then $\Cay(S_n,S)$ is called the {\em modified bubble-sort graph}.
 \\(d) If $T(S)$ is the complete graph $K_n$, then $\Cay(S_n,S)$ is called the {\em complete transposition graph}. 
 \\(e) If $T(S)$ is the complete bipartite graph $K_{k,n-k}$, then $\Cay(S_n,S)$ is called the {\em generalized star graph}. 
\end{Definition}

The automorphism groups of some families of Cayley graphs generated by transpositions have been obtained recently.  We briefly sketch the idea behind the proofs of these results.  In the proofs of these results, one starts with the subgroup $R(S_n) \rtimes \Aut(S_n,S)$ of automorphisms of $X=\Cay(S_n,S)$ and the problem is often to show that this subgroup is in fact the full automorphism group of $X$.   Let $L_e = L_e(X)$ denote the set of automorphisms of $X$ that fixes the vertex $e$ and each of its neighbors.  It can be shown that an equivalent problem is to show that $L_e=1$.   In order to prove $L_e=1$, it suffices to show that an automorphism $g$ of $X$ that fixes the vertex $e$ and each vertex in $X_1(e)$ also fixes each vertex in $X_2(e)$.  For it follows by induction that $g$ fixes all the remaining vertices of $X$.  

One way to prove that an automorphism of $X$ which fixes $e$ and $X_1(e)$ pointwise also fixes $X_2(e)$ pointwise is to establish that there are unique cycles in the graph $X$ that contain the vertex $e$ and some of its neighbors.  Since these cycles are unique, it would then follow that an automorphism of $X$ that fixes $e$ and each of its neighbors must also fix the remaining vertices of the cycle.  

\begin{Lemma} \label{lem:Ganesan:uniqueness:4cycles} \cite{Ganesan:JACO} 
 Let $S$ be a set of transpositions generating $S_n$.  Let $t, k \in S, t \ne k$.  Then, $tk=kt$ if and only if there is a unique 4-cycle in $\Cay(S_n,S)$ containing $e, t$ and $k$. 
\end{Lemma}

Suppose an automorphism $g$ of the Cayley graph $X=\Cay(S_n,S)$ fixes the identity vertex $e$ and each of its neighbors. In other words, suppose $g \in L_e$. Let $t, k \in S$.  If $t$ and $k$ have disjoint support, then by Lemma~\ref{lem:Ganesan:uniqueness:4cycles}, $g$ also fixes all elements of the form $tk$ in $X_2(e)$.   If $t$ and $k$ have overlapping support, then one can often prove that there is a unique 6-cycle in $X$ containing $e, t, k$ and a vertex at distance 3 from $e$. This would imply that $g$ fixes all vertices $tk$ in $X_2(e)$ if $t$ and $k$ have overlapping support.  Hence, $g$ fixes $X_2(e)$ pointwise.  This method is often used to prove that $L_e=1$ and that $R(S_n) \rtimes \Aut(S_n,S)$ contains all automorphisms of $X$.  

Godsil and Royle \cite{Godsil:Royle:2001} proved that a particular family of Cayley graphs of the symmetric group is a family of GRRs:

\begin{Theorem} \label{thm:GR:autgroup:asym:tree}
  Let $S$ be a set of transpositions generating $S_n$ such that the transposition graph $T(S)$ is an asymmetric tree.  Then, the automorphism group of the Cayley graph $\Cay(S_n,S)$ is isomorphic to $S_n$.
\end{Theorem}

Zhang and Huang \cite{Zhang:Huang:2005} obtained the automorphism group of the bubble-sort and modified bubble-sort graphs.  Feng \cite{Feng:2006} generalized Theorem~\ref{thm:GR:autgroup:asym:tree} as follows. 

\begin{Theorem} \label{thm:Feng:tree}
 Let $S$ be a set of transpositions generating $S_n$ such that the transposition graph $T(S)$ is a tree.  Then, the automorphism group of the Cayley graph $\Cay(S_n,S)$ is isomorphic to $R(S_n) \rtimes \Aut(S_n,S)$. 
\end{Theorem}

Feng \cite{Feng:2006} also proved the following result, which holds for arbitrary graphs (not just for trees):

\begin{Theorem}
 Let $S$ be a set of transpositions in $S_n$ $(n \ge 3$).  Then, the set $\Aut(S_n,S)$ of automorphisms of $S_n$ that fixes $S$ setwise is isomorphic to the automorphism group of the transposition graph $T(S)$. 
\end{Theorem}

Ganesan \cite{Ganesan:DM:2013} generalized Theorem~\ref{thm:Feng:tree} from transposition graphs that are trees to arbitrary connected transposition graphs that have girth at least 5:

\begin{Theorem} \label{thm:Ganesan:girth:5}
 Let $S$ be a set of transpositions generating $S_n$ ($n \ge 3$).  If the girth of the transposition graph $T(S)$ is at least 5, then the automorphism group of the Cayley graph $\Cay(S_n,S)$ is isomorphic to $R(S_n) \rtimes \Aut(S_n,S)$. 
\end{Theorem}

Ganesan \cite{Ganesan:AJC:2016} showed that the semidirect product in the expression $R(S_n) \rtimes \Aut(S_n,S)$ for the automorphism group of a normal Cayley graph $\Cay(S_n,S)$ can be strengthened to a direct product.  Since $R(S_n) \cong S_n$ and $\Aut(S_n,S) \cong \Aut(T(S))$, we then obtain

\begin{Theorem}
Let $S$ be a set of transpositions generating $S_n$ such that the Cayley graph $\Cay(S_n,S)$ is normal.  Then, the automorphism group of $\Cay(S_n,S)$ is isomorphic to $S_n \times  
\Aut(T(S))$.
\end{Theorem}

The Cayley graphs investigated in Theorem~\ref{thm:GR:autgroup:asym:tree}, Theorem~\ref{thm:Feng:tree} and Theorem~\ref{thm:Ganesan:girth:5} are normal, i.e. their full automorphism group is equal to $R(S_n) \rtimes \Aut(S_n,S)$.  Ganesan \cite{Ganesan:DM:2013} showed that if the transposition graph is the 4-cycle graph, then the Cayley graph $\Cay(S_n,S)$ is non-normal:

\begin{Theorem} \label{thm:4cycleTS:nonnormal}
 Let $S$ be a set of 4 cyclically adjacent transpositions in $S_4$.  Then, the Cayley graph $\Cay(S_4,S)$ is non-normal. 
\end{Theorem}

When the transposition graph of $S$ is a complete graph, the Cayley graph $\Cay(S_n,S)$ is called the {\em complete transposition graph}.   Ganesan \cite{Ganesan:JACO} determined an equivalent condition for normality and proved that the complete transposition graph is a non-normal Cayley graph.

\begin{Theorem} \label{thm:completeTS:nonnormal}
 Let $S$ be a set of transpositions generating $S_n$. Let $L_e$ be the set of automorphisms of the Cayley graph $X=\Cay(S_n,S)$ that fixes the identity vertex $e$ and each of its neighbors.  Then, the Cayley graph $\Cay(S_n,S)$ is normal iff $L_e=1$.  The complete transposition graph is a non-normal Cayley graph.
\end{Theorem}

The automorphism group of the complete transposition graph was obtained in \cite{Ganesan:JACO}, \cite{Levenshtein:Siemons:2009}:

\begin{Theorem} \label{thm:completeTS:autgroup}
 Let $S$ be the set of all transpositions in $S_n$.  Then, the automorphism group of the Cayley graph $\Cay(S_n,S)$ is equal to 
 \[ \Aut(\Cay(S_n,S))  = (R(S_n) \rtimes \Inn(S_n)) \rtimes \mathbb{Z}_2,\] 
 where $R(S_n)$ is the right regular representation of $S_n$, $\Inn(S_n)$ is the inner automorphism group of $S_n$, and $\mathbb{Z}_2 = \langle h \rangle$, where $h$ is the inverse map $\alpha \mapsto \alpha^{-1}$. 
\end{Theorem}

An open problem is to obtain the automorphism group of the remaining families of Cayley graphs generated by transpositions:  

\begin{Problem}
 Obtain the automorphism group of the remaining families of Cayley graphs generated by transpositions; in particular,  generalize Theorem~\ref{thm:Ganesan:girth:5}. 
\end{Problem}

The only non-normal Cayley graphs generated by transpositions known so far are the Cayley graph generated by the 4-cycle transposition graph (cf. Theorem~\ref{thm:4cycleTS:nonnormal}) and the complete transposition graphs (cf. Theorem~\ref{thm:completeTS:nonnormal}).  We conjecture there are no other non-normal Cayley graphs generated by transpositions: 

\begin{Conjecture}
 Let $S$ be a set of transpositions generating $S_n$ ($n \ge 3$).   If the transposition graph $T(S)$ is not the 4-cycle graph and not the complete graph, then the automorphism group of the Cayley graph $\Cay(S_n,S)$ is isomorphic to $R(S_n) \rtimes \Aut(S_n,S)$. 
\end{Conjecture}

Recently, Ganesan \cite{Ganesan:DMGT:toappear} characterized the isomorphism and edge-transitivity of Cayley graphs generated by transpositions:

\begin{Theorem} 
Let each of $S, S'$ be transposition sets which generate $S_n$.  
The Cayley graphs $\Cay(S_n,S)$ and $\Cay(S_n,S')$ are isomorphic if and only if the transpositions graphs $T(S)$ and $T(S')$ are isomorphic.
\end{Theorem}

\begin{Theorem}  \label{thm:cay:transp:edgetrans}
Let $S$ be a set of transpositions generating $S_n$.  The Cayley graph $\Cay(S_n,S)$ is edge-transitive if and only if the transposition graph $T(S)$ is edge-transitive. 
\end{Theorem}

Recall that if a graph is edge-transitive, then its vertex-connectivity is maximum possible.  Thus, we have the following consequence of Theorem~\ref{thm:cay:transp:edgetrans}:

\begin{Corollary} \label{cor:edgetranscay:kappa:eq:delta}
 Let $S$ be a set of transpositions generating $S_n$.  If the transposition graph $T(S)$ is edge-transitive, then the vertex-connectivity of the Cayley graph $\Cay(S_n,S)$ is maximum possible. In particular, the star graphs, modified bubble-sort graphs and complete transposition graphs have optimal fault-tolerance.
\end{Corollary}

All Cayley graphs generated by transpositions are bipartite, hence are $K_4$-free.  By Theorem~\ref{thm:Mader:K4:free}, all Cayley graphs generated by transpositions have optimal fault-tolerance. Thus, we obtain the following result, which is a generalization of Corollary~\ref{cor:edgetranscay:kappa:eq:delta}. 

\begin{Corollary}
Let $S$ be a set of transpositions generating $S_n$. Then, the vertex-connectivity of $\Cay(S_n,S)$ is maximum possible.
\end{Corollary}

\subsection{Optimal fault-tolerance of some families of Cayley graphs}

In this section, we prove that some families of Cayley graphs have optimal fault-tolerance.  Recall that if $W$ is a minimum separating set for a graph $X$, then the components of $X-W$ are called parts of $X$, and an atomic part of $X$ is a part having the smallest possible number of vertices.  The main result of this section is Theorem \ref{thm:Gao:Novick:atom:in:SS} which gives a necessary condition on the structure of an atomic part of a Cayley graph.  Theorem~\ref{thm:Gao:Novick:atom:in:SS} is then used to prove the optimal fault-tolerance of several families of Cayley graphs, including Cayley graphs generated by transpositions, torus networks, and the folded hypercubes.  

Recall that an atomic part of $X$ is a subgraph of $X$. We shall refer to the vertex set of an atomic part as an {\em atom}.

\begin{Theorem} \cite{Gao:Novick:2007} \label{thm:Gao:Novick:atom:Cayley:is:right:coset}
 Let $H$ be a group and suppose $1 \notin S = S^{-1}$.  Let $A$ be an atom of the Cayley graph $X = \Cay(H,S)$ containing the identity vertex $e$.  Then,
 \\(a) $A$ is a subgroup of $H$.
 \\(b) $A$ is generated by $S \cap A$.
 \\(c) Every atom of $X$ is a right coset of $A$.
\end{Theorem}

\noindent \emph{Proof:}  We first prove (c).

Proof of (c):  Let $r_h \in R(H)$ be the map $(H \rightarrow H, x \mapsto xh)$.  A translate of the atom $A$ is defined to be the image $Ah$ of $A$ under the action of an element $r_h \in R(H)$.  The right regular representation $R(H)$ acts on the vertices of Cayley graph as a group of automorphisms. Hence, for each $h \in H$, the translate $Ah$ of the atom $A$ is also an atom of $X$.  Recall that distinct atoms are disjoint, and the set of atoms of $X$ forms a complete block system for $\Aut(X)$ (cf. Corollary~\ref{cor:atoms:complete:block:system:AutX}).  Thus, the translates of the atom $A$ are the only atoms of $X$.

Proof of (a): Since $H$ is finite, it suffices to show that $A$ satisfies the closure property.  Let $a, b \in A$.  Since $Ab$ and $A$ are atoms of $X$ containing the common element $b$, $Ab=A$.  Hence, $ab \in A$. 

Proof of (b): We show that $A = \langle S \cap A \rangle$.  Suppose $x, y \in A$ are adjacent vertices of $X$, and suppose $y=sx$ for some $s \in S$.  Then, $s=yx^{-1} \in A$ because $A$ is a subgroup. Thus, $s \in A$.  Since $s \in S$, $s \in S \cap A$.  Thus, the generator element (or edge label) corresponding to any two adjacent vertices in $X$ is in $S \cap A$.

The induced subgraph $X[A]$ is connected.  Recall also that $e \in A$.  If $a \in A$, then there is a sequence of edges in $X[A]$, each corresponding to a generator element in $S \cap A$, joining $e$ to $a$.  Hence, $A \le \langle S \cap A \rangle$.  Since $A$ is a subgroup, $A = \langle A \rangle$, and so $A \supseteq \langle S \cap A \rangle$.  Thus, $A = \langle S \cap A \rangle$.
\qed

In what follows, we use the following notation. If $A, B$ are subsets of a group $H$, then $AB := \{ab: a \in A, b \in B\}$.

\begin{Lemma} \label{lem:Gao:Novick:A:B}
 Let $A$ be a group and let $B$ be a subset of $A$.  If $|B| > |A|/2$, then $A = B B^{-1}$.  
\end{Lemma}

\noindent \emph{Proof:}
Since $A$ is a group, by the closure property, $A \supseteq BB^{-1}$.  We need to show the reverse inclusion. Let $a \in A$.  By hypothesis, $|aB| + |B| > |A|$.  Hence, the subsets $aB$ and $B$ have a nontrivial intersection; say $ab_2=b_1$ for some $b_1, b_2 \in B$.  Then, $a = b_1 b_2^{-1} \in BB^{-1}$.  Thus, $A \subseteq B B^{-1}$. 
\qed

\begin{Theorem}
\label{thm:Gao:Novick:atom:in:SS}   \cite{Gao:Novick:2007}
Let $H$ be a group, and suppose $1 \notin S = S^{-1}$.  Let $A$ be the atom of the Cayley graph $X=\Cay(H,S)$ containing the identity vertex $e$.  Then $A \subseteq SS$. 
\end{Theorem}

\noindent \emph{Proof:}
If $\kappa(X) = \delta(X)$, then $|A|=1$, and so $A = \{e\}$. In this case, $A \subseteq S S^{-1} = SS$ and we are done.  For the rest of the proof, suppose $\kappa(X) < \delta(X)$; this implies that $|A| \ge 2$.  We know $A$ is a subgroup of $H$.  Let $N(A)$ denote the set of vertices in $V(X)-A$ that are adjacent to some vertex in $A$.  By definition of $E(X)$, vertices in $sA$ ($s \in S$) are adjacent to vertices in $A$. Hence, $N(A)$ is a union of left cosets $b_1 A, \ldots, b_\ell A$ for some $b_1,\ldots,b_l \in S-A$.  Thus, $\kappa(X) = \ell |A|$.  Note that $S \subseteq A \cup N(A)$ because every vertex in $S$ is adjacent to $e \in A$.  Define $d_0 := |S \cap A|$, $d_i:=|S \cap b_i A|$ for $i=1,\ldots,\ell$.   Then $|S| = d_0+d_1\ldots+d_\ell$.  

We claim that there exists an $i \in \{0,1,\ldots,\ell\}$ such that $d_i > |A|/2$.  Suppose $d_i \le |A|/2$ for $i=1,\ldots,\ell$.  We show that $d_0 > |A|/2$.  We have $|S| = d_0 + d_1 + \ldots + d_\ell \le d_0 + \ell~|A|/2$.  Also, $\ell |A| = \kappa(X) < \delta(X) = |S| \le d_0 + \ell |A|/2$. Hence, $d_0 > \ell |A| - \ell |A|/2 = \ell |A|/2$. This proves the claim.

Let $b_0:=e$.  By the claim in the previous paragraph, $|S \cap b_i A| > |A|/2$ for some $i \in \{0, 1, \ldots, \ell \}$.  Given such an $i$, define the subset $B \subseteq A$ to be such that $S \cap A b_i = B b_i$.  Then, $|B| = |S \cap A b_i| > |A|/2$.  By Lemma~\ref{lem:Gao:Novick:A:B}, $A = B B^{-1} = B b_i b_i^{-1} B^{-1} = B b_i (B b_i)^{-1} \subseteq S S^{-1} = SS$.
\qed

Theorem~\ref{thm:Gao:Novick:atom:in:SS} was used by Gao and Novick \cite{Gao:Novick:2007} to prove that all Cayley graphs of the symmetric group generated by transpositions have optimal fault-tolerance.

\begin{Corollary} 
 Let $S$ be a set of transpositions generating $S_n$. Then, the Cayley graph $X=\Cay(S_n,S)$ has optimal fault-tolerance. 
\end{Corollary}

\noindent \emph{Proof:}
By way of contradiction, suppose $\kappa(X) < \delta(X)$.  Then, there exists an atom $A$ such that $e \in A$ and $|A| \ge 2$.  By Theorem~\ref{thm:Gao:Novick:atom:in:SS}, $A \subseteq SS$.  By Theorem~\ref{thm:Gao:Novick:atom:Cayley:is:right:coset}, $A = \langle A \cap S \rangle$. Since $|A| \ge 2$, $A \cap S$ is nonempty.  It follows that some generating transposition $t \in S$ is in $A$, and hence in $S S$.  This implies that $S$ and $SS$ have a nontrivial intersection, which is impossible because $S$ contains a set of odd permutations and $SS$ a set of even permutations. 
\qed

\begin{Exercise}
 Let $X$ be the 3-dimensinal torus graph $\Cay( \mathbb{Z}_r \times \mathbb{Z}_s \times \mathbb{Z}_t, S)$, where $S = \{ \pm (1,0,0), \pm(0,1,0), \pm(0,0,1)\}$.  Observe that $X$ is a 6-regular graph having $rst$ vertices.  Prove that the vertex-connectivity of $X$ is equal to its valency. 
\end{Exercise}

\subsection{Cayley graphs from linear codes}

In this section, we introduce Cayley graphs generated by matrices, especially parity-check matrices of linear codes.  We show how many families of Cayley graphs, including the hypercubes and folded hypercubes, arise as  special cases of Cayley graphs associated with matrices, and we prove the optimal fault-tolerance of some families of Cayley graphs.  Let $\mathbb{Z}_2^r$ denote the abelian group which consists of the set of all 0-1 vectors of length $n$, with the group operation being vector addition (thus, $\mathbb{Z}_2^r$ is the abelian group associated with the vector space $\mathbb{F}_2^r$).  We first apply Theorem~\ref{thm:Gao:Novick:atom:in:SS} to obtain a sufficient condition for a Cayley graph of the group $\mathbb{Z}_2^r$ to have optimal connectivity.   

Given any $r \times n$ matrix $\mathcal{H}$ whose entries are over the binary field $\mathbb{F}_2$, we can form a set $S$ consisting of the columns of $\mathcal{H}$.  The {\em  Cayley graph generated by the matrix $\mathcal{H}$} is defined to be the Cayley graph $\Cay(\mathbb{Z}_2^r,S)$.  Conversely, given any subset $S \subseteq \mathbb{Z}_2^r$, we can form a matrix $\mathcal{H}$ whose columns are the elements of $S$ (in any order).  Thus, every Cayley graph of the group $\mathbb{Z}_2^r$ is isomorphic to the Cayley graph generated by some 0-1 matrix $\mathcal{H}$, and conversely. 

In order to investigate the applications of coding theory to interconnection networks, we shall especially consider matrices $\mathcal{H}$ which are parity-check matrices of linear codes. For this reason, we assume the matrices $\mathcal{H}$ are short and fat matrices, of dimension $r \times n$ where $r \le n$.  We also assume that the matrix $\mathcal{H}$ has rank $r$; in this case, if $S$ is the set of columns of $\mathcal{H}$, then the Cayley graph $\Cay(\mathbb{Z}_2^r,S)$ is a connected graph because the columns of $\mathcal{H}$ span the entire $r$-dimensinal vector space.  We state the above remarks formally:

\begin{Definition}
 Let $\mathcal{H}$ be an $r \times n$ matrix with entries over the binary field $\mathbb{F}_2$, where $r \le n$ and $rank(H)=r$. Let $S$ be the set of columns of $\mathcal{H}$.  The {\em Cayley graph generated by $\mathcal{H}$} is defined to be the Cayley graph $\Cay(\mathbb{Z}_2^r,S)$. 
\end{Definition}

We assume that all columns of $\mathcal{H}$ are nonzero; this is equivalent to assuming that the generating set $S$ of the Cayley graph does not contain the identity element, which implies that the Cayley graph has no self-loops. Also, it is clear that each element of $S$ is its own inverse since we are working over the binary field.  Thus, the Cayley graph generated by a matrix $\mathcal{H}$ is a simple, undirected graph.  By definition of the Cayley graph of a matrix $\mathcal{H}$,  the order in which the columns of $\mathcal{H}$ appear is irrelevant.  

In the special case where $\mathcal{H}$ is the $r \times r$ identity matrix, the Cayley graph generated by $\mathcal{H}$ is isomorphic to the hypercube $Q_r$.  In fact, the Cayley graph generated by any square matrix having full rank is isomorphic to the hypercube:

\begin{Lemma}
Let $\mathcal{H}$ be an $r \times r$ 0-1 matrix having rank $r$. Then, the Cayley graph generated by $\mathcal{H}$ is isomorphic to the hypercube graph $Q_r$. 
\end{Lemma}

\noindent \emph{Proof:}
Let $X$ be the Cayley graph generated by $\mathcal{H}$. Let $\theta$ be the linear transformation from $\mathbb{F}_2^n$ to itself defined by $\theta: x \mapsto \mathcal{H}x$ .  Since $\mathcal{H}$ has full rank, $\theta$ is a bijection from $V(Q_n)$ to $V(X)$.  Also, $\theta$ is an isomorphism from the graph $Q_n$ to the graph $X$ because it preserves adjacency and non-adjacency: $u + v = e_i$ for some $i$ if and only if $\mathcal{H}(u+v)$ is a column of $\mathcal{H}$.
\qed

More generally \cite{Godsil:Royle:2001},

\begin{Exercise}
If $H$ is a group, $S \subseteq H$ and $\theta \in \Aut(H)$, then the Cayley graphs $\Cay(H,S)$ and $\Cay(H, \theta(S))$ are isomorphic. 
\end{Exercise}

\begin{Example}
In the special case where $\mathcal{H}$ is the $4 \times 5$ matrix
\[
 \mathcal{H} = \left[ \begin{array}{ccccc}
                       1 & 0 & 0 & 0 & 1 \\
                       0 & 1 & 0 & 0 & 1 \\
                       0 & 0 & 1 & 0 & 1 \\
                       0 & 0 & 0 & 1 & 1 \\
                      \end{array}
\right],
\]
the Cayley graph generated by $\mathcal{H}$ is the folded hypercube graph $FQ_4$.  Observe that $\mathcal{H}$ is the parity-check matrix of the repetition code of dimension 1 and length 5 (cf. \cite{MacWilliams:Sloane:1977}).  If $\mathcal{H}$ is the $r \times n$ parity-check matrix of the repetition code of length $n$ and dimension $1$ (so $r=n-1$), then the Cayley graph generated by $\mathcal{H}$ is isomorphic to the folded hypercube $FQ_r$. 
\qed
\end{Example}

More generally, if $\mathcal{H}$ is an $r \times n$ matrix with $r < n$ and having rank $r$, then the Cayley graph generated by $\mathcal{H}$ is some supergraph of the hypercube $Q_r$.   

We now give a sufficient condition for optimal vertex-connectivity of Cayley graphs of $\mathbb{Z}_2^r$.

\begin{Corollary} \label{cor:suff:cond:matrixH:optfaulttol}
 Let $\mathcal{H}$ be an $r \times n$ matrix of rank $r$ over the binary field.  Let $S$ be the set of columns of $\mathcal{H}$. If no column of $\mathcal{H}$ is a sum of any two columns of $H$, then the Cayley graph $\Cay(\mathbb{Z}_2^r,S)$ has optimal fault-tolerance.
\end{Corollary}

\noindent \emph{Proof:}
Let $X = \Cay(\mathbb{Z}_2^r,S)$.  Suppose $\kappa(X) < \delta(X)$.  Then, there exists an atom $A$ containing the identity vertex $e$ such that $|A| \ge 2$.  We know $A = \langle A \cap S \rangle$.  So $A \cap S$ is nonempty.  Also, $A \subseteq S+S := \{x+y: x, y \in S\}$.  Hence, some element of $S$ must be in $S+S$. But this implies that $S$ and $S+S$ overlap, i.e. that some column of $\mathcal{H}$ is a sum of two columns of $\mathcal{H}$.  
\qed

In addition to the proofs given earlier of the optimal fault-tolerance of the hypercubes and folded hypercubes, we now give another proof that the hypercube and the folded hypercube graphs have optimal fault-tolerance.

\begin{Corollary}
 The hypercube graphs $Q_n$ and folded hypercube graphs $FQ_n$ ($n \ge 4$) have optimal fault-tolerance. 
\end{Corollary}

\noindent \emph{Proof:}
Observe that the $r \times n$ matrix $\mathcal{H}$ generating these two families of Cayley graphs satisfies the condition that no column is a sum of two columns.  More specifically, define the weight of a 0-1 vector to be the number of 1's in the vector.  For the hypercube $Q_n$, the matrix $\mathcal{H}$ is the identity matrix.  Every column of $\mathcal{H}$ is a vector of weight 1, while the sum of any two columns of $\mathcal{H}$ has weight 0 or 2.  By Corollary~\ref{cor:suff:cond:matrixH:optfaulttol}, $Q_n$ has optimal fault-tolerance.  The proof for $FQ_n$ is left to the reader. 
\qed

%
A general research direction is to investigate the properties of Cayley graphs generated by $\mathcal{H}$ for various families of matrices $\mathcal{H}$.  Some results in this direction where $\mathcal{H}$ is a parity-check matrix or generator matrix of a linear code are in 
\cite{Zemor:Cohen:1992}
\cite{Bruck:Cypher:Ho:1993}
\cite{Bruck:Ho:1996}
\cite{Bruck:Cypher:Ho:1995}
and \cite{Garcia:Peyrat:1997}.
%
%
%
{
\bibliographystyle{unsrt}
\bibliography{refsaut}
}

\end{document}